# ESTIMATING THE DEGREE OF ACTIVITY OF JUMPS IN HIGH FREQUENCY DATA

By Yacine Aït-Sahalia[1] and Jean Jacod

*Princeton University and UPMC (Université Paris-6)*

We define a generalized index of jump activity, propose estimators of that index for a discretely sampled process and derive the estimators' properties. These estimators are applicable despite the presence of Brownian volatility in the process, which makes it more challenging to infer the characteristics of the small, infinite activity jumps. When the method is applied to high frequency stock returns, we find evidence of infinitely active jumps in the data and estimate their index of activity.

**1. Introduction.** Using high frequency financial data, which are now widely available, we can hope to answer a number of questions regarding the characteristics of the process that drives asset returns. Let us model the log-price $X$ of some asset as a 1-dimensional process, which we will observe over a fixed time interval $[0, T]$ at discrete times $0, \Delta_n, 2\Delta_n, \ldots$ with a time interval $\Delta_n$ between successive observations that is small. This is the essence of high frequency data. Let us further assume that this process is an Itô semimartingale, meaning that its characteristics are absolutely continuous with respect to Lebesgue measure. So, it has a drift, a continuous martingale part that is the integral of a possibly stochastic process with respect to a Brownian motion, and we will also let it have jumps with a possibly stochastic Lévy measure.

For modeling purposes, one would like to infer the characteristics of $X$ from the observations; that is, its drift, volatility and Lévy measure. When the time interval $\Delta_n$ goes to 0, it is well known that one can consistently infer the volatility under very weak assumptions. However, such consistent inference is impossible for the drift or the Lévy measure, if the overall time interval $[0, T]$ is kept fixed.

Received March 2008; revised July 2008.
[1]Supported in part by NSF Grant DMS-05-32370.
*AMS 2000 subject classifications.* Primary 62F12, 62M05; secondary 60H10, 60J60.
*Key words and phrases.* Jumps, index of activity, infinite activity, discrete sampling, high frequency.

In fact, even in the unrealistic case where the whole path of $X$ is observed over $[0, T]$, one can infer neither the drift nor the Lévy measure. One can, however, hope to be able to characterize the behavior of the Lévy measure near 0: first, whether it does not explode near 0, meaning that the number of jumps is finite; and second, when this number is infinite, we would like to be able to say something about the concentration of small jumps. Our objective in doing so is to provide specification tools for financial models, where the presence or at least possibility of large jumps is generally accepted. There is much less consensus in the literature regarding the nature or even the need for small jumps.

For this purpose, let us define, for a generic semimartingale $X$,

$$(1) \quad B(r)_t = \sum_{s \leq t} |\Delta X_s|^r, \qquad I_t = \{r \geq 0 : B(r)_t < \infty\}, \qquad \beta_t = \inf(I_t),$$

where $\Delta X_s = X_s - X_{s-}$ is the size of the jump at time $s$, and $r \geq 0$, with the convention $0^0 = 0$. Necessarily, the (random) set $I_t$ contains the interval $(\beta_t, \infty)$, whereas it may contain $\beta_t$ itself or not. Moreover, $2 \in I_t$ always, and, of course, $t \mapsto \beta_t$ is nondecreasing. Hence, if we observe the whole path of $X$ over $[0, T]$, we know the sets $I_t(\omega)$ and the numbers $\beta_t(\omega)$ for all $t \leq T$.

We call $\beta_T(\omega)$ the *jump activity index* for the path $t \mapsto X_t(\omega)$ at time $T$ (or, more precisely, up to time $T$). We define this index in analogy with the special case where $X$ is a Lévy process. In this case, $I_t$ and $\beta_t$ are no longer random. Further, they do not depend on the time $t$, and $I_t$ is also the set of all $r \geq 0$ such that $\int_{\{|x| \leq 1\}} |x|^r F(dx) < \infty$, where $F$ is the Lévy measure. This property shows that, for a Lévy process, the jump activity index coincides with the Blumenthal–Getoor index of the process [see Blumenthal and Getoor (1961)]. In the further special case where $X$ is a stable process, $\beta$ is also the stable index of the process.

When $X$ is a Lévy process, the interval $I$ and the index $\beta$ are, of course, only tiny elements of the whole Lévy measure $F$, which convey approximately the same information ($I$ gives slightly more information than $\beta$). However, the value of $\beta$ is probably the most informative knowledge one can draw about $F$ from the observation of the path $t \mapsto X_t$ for all $t \leq T$ when $T$ is finite. Things are very different when $T \to \infty$, though, since observing $X$ over $[0, \infty)$ completely specifies $F$. But, when the time horizon $T$ is kept fixed, and with the whole path observed over $[0, T]$, we can infer only the behavior of the Lévy measure $F$ near 0 (because we need a potentially infinite number of observations for consistent estimation). Then, $\beta$ captures an essential qualitative feature of $F$, which is its level of activity, which is that when $\beta$ increases, the (small) jumps tend to become more and more frequent.

$\beta$ is related to the "degree of activity" of jumps. All Lévy measures put finite mass on the set $(-\infty, -\varepsilon] \cup [\varepsilon, +\infty)$ for any arbitrary $\varepsilon > 0$; therefore, if the process has infinite jump activity, then it must be because of the "small" jumps, which are defined as those smaller than $\varepsilon$. If $F([-\varepsilon,\varepsilon]) < \infty$, then the process has finite activity and $0 \in I$, or, equivalently, $\beta = 0$. But, if $F([-\varepsilon,\varepsilon]) = \infty$, then the process has infinite activity, and, in addition, $\beta > 0$ as long as the Lévy measure $F([-\varepsilon,\varepsilon])$ diverges near 0 at a rate larger than a power $\varepsilon^{-a}$ for some $a > 0$. The higher $\beta$ gets (up to 2), the more active the small jumps become. The same remarks also apply for general semimartingales. These properties are what motivate our calling $\beta$ a jump activity index and our interest in estimating it.

In the more realistic situation where the semimartingale $X$ is only observed at times $i\Delta_n$ over $[0,T]$, the estimation problem is made more challenging by the presence in $X$ of a continuous martingale part. By its very nature, $\beta_T$ characterizes the behavior of $F$ near 0. Hence, it is natural to expect that the small increments of the process are going to be the ones that are most informative about $\beta_T$. But, those small increments are precisely the ones where the contribution from the continuous martingale part of the process is inexorably mixed with the contribution from the small jumps. Being able to "see through" the continuous part of the semimartingale in order to say something about the number and concentration of small jumps is going to be the challenge we face as we attempt to estimate $\beta_T$.

Related to this paper are Woerner (2006), who proposes an estimator of the jump activity index and of the Hurst exponent, but in the absence of a continuous Brownian part to the semimartingale, and Cont and Mancini (2007), who propose a test for the finiteness of the variation of the jump part. Also, Belomestny (2008) estimates the same index when there is no Brownian part and when, together with the prices, some option prices also are recorded. Another related problem is the estimation of the index $\beta$ of a stable process [see, e.g., DuMouchel (1983)]. However, our situation here is fundamentally different from those, in that we also have a continuous part in the semimartingale. The situation is also different from that in Aït-Sahalia and Jacod (2008), where we studied Fisher's information for the parameters of a Brownian plus stable pure jump process, but the jump process was the dominant component in that paper. Here, the continuous part of the semimartingale dominates the small increments, and we estimate the activity index of a pure jump process where the dominant component is a Brownian motion.

The aim of this paper is to construct estimators $\hat{\beta}_n(T)$ for $\beta_T$, which are consistent when $\Delta_n \to 0$, and to provide rates of convergence and asymptotic distributions. Ideally, we would also like to have estimators that are, as much as possible, model-free, in the sense that they behave well without too strong assumptions on the form of the drift, the volatility or the Lévy measure.

As it turns out, a fully model-free behavior of the estimators may be too much to ask. The assumptions we make below on the drift and the volatility process are quite unrestrictive, but obtaining rates of convergence will require more specific assumptions on the Lévy measure. In particular, we will assume that the main part of the Lévy measure near 0 behaves locally like the Lévy measure of a stable process, and we will provide estimators and their properties when the index is $\beta > 0$. This assumption seems to be unavoidable, since, as we shall also see, even when $X$ is a Lévy process, strong assumptions on the Lévy measure are necessary. At this juncture, it may be worth noting that considering semimartingales rather than simply Lévy processes, or exponentials of such, does not change or weaken the results.

The paper is organized as follows. In Section 2, we formally define the index of jump activity, construct estimators for it and present the main properties of the estimators in the general case where the process is a semimartingale. Section 3 is devoted to the special and simpler case of a symmetric stable process, and Section 4 is about more general Lévy processes. We propose a small sample bias correction in Section 5. We present the results of Monte Carlo simulations in Section 6, and we compute our estimators over all 2006 transactions of the Dow Jones stocks in Section 7, focusing in particular on Intel and Microsoft. Section 8 is devoted to technical results and to the proof of the main theorems, which apply to Itô semimartingales, under suitable assumptions on the Lévy measure.

## 2. The model and main results.

2.1. *Defining an index of jump activity.* Our structural assumption is that $X$ is a 1-dimensional Itô semimartingale on some filtered space $(\Omega, \mathcal{F}, (\mathcal{F}_t)_{t \geq 0}, \mathbb{P})$, which means that its characteristics $(B, C, \nu)$ are absolutely continuous with respect to Lebesgue measure [see Jacod and Shiryaev (2003) for all notions not explained here]. In other words, the characteristics of $X$ have the form

$$(2) \qquad B_t = \int_0^t b_s \, ds, \qquad C_t = \int_0^t \sigma_s^2 \, ds, \qquad \nu(dt, dx) = dt \, F_t(dx).$$

Here, $b = (b_t)$ and $\sigma = (\sigma_t)$ are real-valued optional processes, and $F_t = F_t(\omega, dx)$ is a predictable random measure, meaning that for all Borel sets $A$ in $\mathbb{R}$ the process $(F_t(A))$ is predictable (possibly taking the value $+\infty$). This model is quite general. For instance, the drift, volatility and jump measures can be stochastic and jump themselves.

There are other ways of expressing this assumption, for example through a Wiener process $W$ and a Poisson random measure $\underline{\mu}$ with compensator $\underline{\nu}(dt, dx) = dt \otimes dx$ (up to a possible enlargement of the space), as

$$X_t = X_0 + \int_0^t b_s \, ds + \int_0^t \sigma_s \, dW_s$$

$$(3) \qquad + \int_0^t \int_{\mathbb{R}} \delta(s,x) 1_{\{|\delta(s,x)| \leq 1\}} (\underline{\mu} - \underline{\nu})(ds,dx)$$

$$+ \int_0^t \int_{\mathbb{R}} \delta(s,x) 1_{\{|\delta(s,x)| > 1\}} \underline{\mu}(ds,dx).$$

In this formulation, $b$ and $\sigma$ are the same as in (2), $\delta = \delta(\omega,t,x)$ is a predictable function and the connection with $F_t$ is that $F_t(\omega,dx)$ is the restriction to $\mathbb{R} \setminus \{0\}$ of the image of the Lebesgue measure by the map $x \mapsto \delta(\omega,t,x)$. However, it is easier for the problem at hand to express the assumptions on $F_t$ rather than on $\delta$, which, moreover, is not unique (whereas $F_t$ is uniquely defined, up to null sets).

Below, for any measure $H$ on $\mathbb{R}$ we denote by $\overline{H}$ its (symmetrical) tail function

$$(4) \qquad x > 0 \mapsto \overline{H}(x) = H([-x,x]^c).$$

Observe that $B(r)_t < \infty$ if and only if $B'(r)_t < \infty$, where $B'(r)_t = \sum_{s \leq t} |\Delta X_s|^r \wedge 1$, and the process $B'(r)$ is finite-valued if and only if it is locally integrable. In other words, divergence, when it occurs, is caused by the small jumps. Moreover, for any stopping time $T$, we have

$$\mathbb{E}(B'(r)_T) = \mathbb{E}\left( \int_0^T ds \int_{\mathbb{R}} (|x|^r \wedge 1) F_s(dx) \right).$$

We call *instantaneous jump activity index* at time $t$ the (random) number

$$(5) \qquad \beta_t^i = \inf \left\{ r > 0 \colon \int_{\mathbb{R}} (|x|^r \wedge 1) F_s(dx) < \infty \right\}.$$

In light of (5), this is a natural generalization of the notion of Blumenthal–Getoor index for Lévy processes [see Blumenthal and Getoor (1961)]. $\beta^i$ is a predictable process taking its values in $[0,2]$. This process is also characterized by the property that, for any $\varepsilon > 0$, we have

$$(6) \qquad \lim_{x \to 0} x^{\beta_t^i + \varepsilon} \overline{F}_t(x) = 0, \qquad \limsup_{x \to 0} x^{\beta_t^i - \varepsilon} \overline{F}_t(x) = \infty.$$

The "lim sup" above is usually not a limit.

Note, finally, that $\beta_t^i = 0$ does not necessarily imply that the process has finite jump activity, since it is possible for the Lévy measure to diverge slowly near 0, at a subgeometric speed. An example of this would be the Gamma process, which has $F(dx) = (\eta \exp(-\kappa x) 1_{\{x>0\}}/x) \, dx$, so that the Lévy measure $\overline{F}_t(\varepsilon)$ diverges at a logarithmic rate in $\varepsilon$. Finite activity processes (compound Poisson) will have $\beta_t^i = 0$ a.s., however.

2.2. *Assumptions.* We make two assumptions. The first one, on the drift $b$ and volatility $\sigma$, is quite mild.

ASSUMPTION 1. The processes $b$ and $\sigma$ are locally bounded.

The second assumption, on the Lévy measures $F_t$, is more specific. Essentially, we split $F_t$ as $F_t = F_t' + F_t''$, where:

- $F_t'$ is very close to the Lévy measure of a $\beta$-stable process, restricted to a random interval $(-z_t^{(-)}, z_t^{(+)})$ around 0 with some $\beta$ that is not random; the random interval may be empty for some $(\omega, t)$, but not for all;
- $F_t''$ is another Lévy measure with jump activity index less than some $\beta' < \beta$.

The precise statement of the assumption is as follows.

ASSUMPTION 2. There are three (nonrandom) numbers $\beta \in (0, 2)$, $\beta' \in [0, \beta)$ and $\gamma > 0$, and a locally bounded process $L_t \geq 1$, such that we have, for all $(\omega, t)$,

$$(7) \qquad\qquad F_t = F_t' + F_t'',$$

where:

(a) $F_t'$ has the form

$$(8) \quad F_t'(dx) = \frac{1 + |x|^\gamma f(t, x)}{|x|^{1+\beta}} (a_t^{(+)} 1_{\{0 < x \leq z_t^{(+)}\}} + a_t^{(-)} 1_{\{-z_t^{(-)} \leq x < 0\}})\, dx$$

for some predictable nonnegative processes $a_t^{(+)}, a_t^{(-)}, z_t^{(+)}$ and $z_t^{(-)}$ and some predictable function $f(\omega, t, x)$ satisfying

$$(9) \qquad \begin{cases} \dfrac{1}{L_t} \leq z_t^{(+)} \leq 1, & \dfrac{1}{L_t} \leq z_t^{(-)} \leq 1, \qquad a_t^{(+)} + a_t^{(-)} \leq L_t, \\ 1 + |x| f(t, x) \geq 0, & |f(t, x)| \leq L_t; \end{cases}$$

(b) $F_t''$ is a measure that is singular with respect to $F_t'$ and satisfies

$$(10) \qquad\qquad \int_{\mathbb{R}} (|x|^{\beta'} \wedge 1)\, F_t''(dx) \leq L_t.$$

We will also need the increasing and locally bounded process

$$(11) \qquad \bar{A}_t = \int_0^t A_s\, ds, \qquad \text{where } A_t = \frac{a_t^{(+)} + a_t^{(-)}}{\beta}.$$

REMARK 1. In view of (6), the instantaneous index at time $t$ "due to" part $F_t'$ of the Lévy measure is $\beta$ on the set $\{A_t > 0\}$, and 0, otherwise; whereas, the one due to $F_t''$ is everywhere smaller than $\beta'$. Hence, outside a null set, we have $\beta_t = \beta$ on the set $\{\bar{A}_t > 0\}$ and $\beta_t \leq \beta'$, otherwise.

REMARK 2. One could formulate Assumption 2 slightly differently by writing $F_t = F_t^1 + F_t^2$, where $F_t^1$ is given by (8) with $f(t,x) = 0$ and $F_t^2$ satisfying (10) with some $\beta''$, and, further, the restriction of $F_t^2$ to $[-z^{(-)}, z^{(+)}]$ has an absolutely continuous part with a density of the form $f(t,x)|x|^{\gamma-1-\bar{\beta}}$ with (9). The two formulations are equivalent, provided that we take $\beta' = \beta'' \vee (\beta - \gamma)$.

REMARK 3. Take any process of the form

$$(12) \qquad dX_t = b_t \, dt + \sigma_t \, dW_t + \delta_{t-} \, dY_t + \delta_{t-}' \, dY_t',$$

where $\delta$ and $\delta'$ are càdlàg adapted processes, $Y$ is $\beta$-stable or tempered $\beta$-stable and $Y'$ is any other Lévy process whose Lévy measure integrates $|x|^{\beta'}$ near the origin and has an absolutely continuous part whose density is smaller than $K|x|^{\gamma-1-\beta}$ on $[-1,1]$ for some $\gamma > 0$ (e.g., a stable process with index strictly smaller than $\beta'$). Then, $X$ will satisfy Assumption 2. For instance, when further $Y$ is symmetrical with Lévy density $D/|x|^{1+\beta}$, it is satisfied with the two numbers $\beta$ and $\beta'$, $\gamma$ (as above), $f(t,x) = 0$, $z_t^{(-)} = z_t^{(+)} = 1$ and $a_t^{(-)} = a_t^{(+)} = D\delta_{t-}^\beta$.

REMARK 4. When $X$ is a Lévy process, so that $F_t(\omega, dx) = F(dx)$ does not depend on $t$ and $\omega$, Assumption 2 is related with the property that $X$ is a "regular Lévy process of exponential type," as introduced in Boyarchenko and Levendorskiĭ (2002), with $\beta = \nu$ and $\beta' \wedge (\beta - \gamma) = \nu'$ in the notation of that paper. These two assumptions are not exactly comparable. The one in Boyarchenko and Levendorskiĭ (2002) is more stringent about the behavior of "big" jumps, whereas ours is slightly more demanding for "small" jumps.

2.3. *The estimators.* Recall that we observe $X_{i\Delta_n}$ for $i = 0, 1, \ldots, [T/\Delta_n]$. While the processes $B(r)$ are defined from the jumps $\Delta X_s = X_s - X_{s-}$ of $X$, we do not observe these jumps directly. Rather, all that we observe are the discrete increments

$$(13) \qquad \Delta_i^n X = X_{i\Delta_n} - X_{(i-1)\Delta_n}.$$

From these increments, we could try to evaluate $B(r)_T$ and then infer $\beta$. Finding consistent estimators for $B(r)_T$ is easy, but deducing from them an estimator for $\beta$ is almost impossible, because we need to decide whether $B(r)_T$ is infinite or not based on a finite sample.

So, we propose the following idea. For fixed $\varpi > 0$ and $\alpha > 0$, we write

$$(14) \qquad U(\varpi, \alpha)_t^n = \sum_{i=1}^{[t/\Delta_n]} 1_{\{|\Delta_i^n X| > \alpha \Delta_n^\varpi\}}$$

for the number of increments whose magnitude is greater than $\alpha \Delta_n^\varpi$. In all cases below, we will set $\varpi < 1/2$.

To better understand our rationale for doing this, consider the special case $X = \sigma W + Y$, where $Y$ is a $\beta$-stable process, so $\beta_t(\omega) = \beta$. Any increment $\Delta_i^n X = X_{i\Delta_n} - X_{(i-1)\Delta_n}$ satisfies $\Delta_i^n X = \sigma \Delta_n^{1/2} W_1 + \Delta_n^{1/\beta} Y_1$ (equality in law). Then, recalling that $\beta < 2$ and $\Delta_n \to 0$, with a large probability $\Delta_i^n X$ is close to $\sigma \Delta_n^{1/2} W_1$ in law. Those increments give essentially no information on $Y$ and are of order of magnitude $\Delta_n^{1/2}$. However if $Y$ has a "big" jump at time $s$, the corresponding increment is close to $\Delta Y_s$. Hence, one has to throw away all the "small" increments. However, $\beta$ is related to the behavior of $F$ near 0 and, hence, to the "very small" jumps of $Y$. This is why we will use only increments bigger than a cutoff $\alpha \Delta_n^\varpi$ for some $\varpi \in (0, 1/2)$. Asymptotically, those increments are big, because, since $\Delta_n^{1/2} \ll \Delta_n^\varpi$, the main contribution is due to $Y$. Those increments mostly contain a single "big" jump of size of order at least $\Delta_n^\varpi$, and we still get some information on small jumps, because $\Delta_n^\varpi \to 0$.

So, by using the statistic $U$, which simply counts the number of large increments, defined as those greater than $\alpha \Delta_n^\varpi$, we are retaining only those increments of $X$ that are not predominantly made of contributions from its continuous semimartingale part, which are $O_p(\Delta_n^{1/2})$, and instead are predominantly made of contributions due to a jump. The same heuristics work for more general Itô semimartingales.

As we will see later, the key property of the functionals $U(\varpi, \alpha, \Delta_n)$ is their convergence in probability

$$(15) \qquad \Delta_n^{\varpi \beta} U(\varpi, \alpha)_t^n \overset{\mathbb{P}}{\longrightarrow} \frac{\bar{A}_t}{\alpha^\beta},$$

which we will show holds under Assumption 2. This property leads us to propose an estimator of $\beta$ at each stage $n$. Fix $0 < \alpha < \alpha'$ and define

$$(16) \qquad \hat{\beta}_n(t, \varpi, \alpha, \alpha') = \frac{\log(U(\varpi, \alpha)_t^n / U(\varpi, \alpha')_t^n)}{\log(\alpha'/\alpha)},$$

which is at least consistent for estimating $\beta$ on the set $\{\bar{A}_t > 0\}$. If either value of $U$ in (16)–(17) is 0, then, by convention, we set the estimator to be 0.

$\hat{\beta}_n$ is constructed from a suitably scaled ratio of two $U$'s evaluated on the same time scale $\Delta_n$ at two different fixed levels of truncation of the increments $\alpha$ and $\alpha'$. In a way, this construction is in the same spirit as the classical estimator of Hill (1975), who conducts inference about the tails of a distribution based on ratios of various extremes.

We can also propose a second estimator defined as

$$(17) \qquad \hat{\beta}'_n(t, \varpi, \alpha) = \frac{\log(U(\varpi, \alpha)^n_t / U_2(\varpi, \alpha)^n_t)}{\varpi \log 2},$$

where $U_2(\varpi, \alpha)^n_t$ is defined analogously to $U(\varpi, \alpha)^n_t$ in (14), except that sampling at $\Delta_n$ is replaced by sampling at $2\Delta_n$. That is, $\hat{\beta}'_n$ is constructed from a suitably scaled ratio of two $U$'s evaluated at the same level of truncation $\alpha$ on two separate time scales $\Delta_n$ and $2\Delta_n$.

One could also look at a third estimator $\hat{\beta}''_n$ obtained from two $U$'s evaluated at two different rates of truncation $\varpi$ and $\varpi'$. One could further consider estimators based not just on counting the increments that exceed a certain cutoff but also on the magnitude of these increments, as in the case of power variations truncated to use only the large increments.

In the rest of the paper, we will focus mainly on the properties of the estimator $\hat{\beta}_n$, noting that a similar type of analysis yields the consistency and asymptotic distribution of the other two estimators. In general, the asymptotic variance of $\hat{\beta}_n$ is smaller than that of $\hat{\beta}'_n$ and $\hat{\beta}''_n$.

Before studying the properties of the estimator $\hat{\beta}_n$, let us make a few remarks.

REMARK 5. Asymptotically, as $n \to \infty$, the above estimators behave well. However, for any given $n$, it may happen that they are not informative, because too few increments are retained, up to the extreme case where $U(\varpi, \alpha)^n_t = U(\varpi, \alpha')^n_t = 0$. If this is the case, one should take smaller values of $\alpha$ and $\alpha'$.

REMARK 6. Even when the estimators are well defined, they may take a value bigger than or equal to 2. In this case, the estimation is not reliable, and it may be an indication that Assumption 2 is simply not satisfied, which would be the case for example if there is no jump at all in the observed path. So it would make sense to convince oneself that jumps are present [see, e.g., Aït-Sahalia and Jacod (2009)] before attempting to estimate $\beta$.

REMARK 7. As we will see below, asymptotic considerations lead to the selection of $\varpi = 1/5$ as a universal choice valid for all possible values of $\beta$. The cutoff for large increments is $\alpha \Delta_n^\varpi$. When implementing the estimator in practice, in any given sample the value of $\Delta_n$ is fixed, so $\varpi$ and $\alpha$ are not independent parameters. The level of truncation $\alpha$ may be set in relation to the volatility of the continuous part of the semimartingale [i.e., $(t^{-1} \int_0^t \sigma_s^2 \, ds)^{1/2}$] since the objective is to eliminate the increments that are mainly due to the continuous part. The truncation level can be selected in a data-driven manner. Despite the presence of jumps, that volatility can be estimated using the small increments of the process, since

$$(18) \qquad \sum_{i=1}^{[t/\Delta_n]} |\Delta_i^n X|^2 1_{\{|\Delta_i^n X| \leq \alpha \Delta_n^{\varpi}\}} \xrightarrow{\mathbb{P}} \int_0^t \sigma_s^2 \, ds$$

for any $\alpha > 0$ and $\varpi \in (0, 1/2)$. We can then set the cutoff level $\alpha$ to yield a number of (estimated) standard deviations of the continuous part of the semimartingale. For the estimator $\hat{\beta}_n$, $\alpha'$ can then be set as a multiple of $\alpha$. These data-driven choices determine a range of reasonable values for $(\alpha, \alpha')$. One possibility is then to simply average the estimators $\hat{\beta}_n$ obtained for the values of $(\alpha, \alpha')$ over that range. The parameters $(\alpha, \alpha')$ effectively play a role similar to that of bandwidth parameters in a nonparametric analysis.

REMARK 8. The construction of the estimators relies on the property (15), which holds under (slightly) weaker assumptions than Assumption 2, provided that the definition of $\bar{A}_t$ and the rate of convergence are suitably amended. As a result, the estimators given in (16) remain consistent under weaker assumptions. For example, when $X$ has only finitely many jumps, the index is $\beta = 0$, and $U(\varpi, \alpha)_t^n$ converges to the number of jumps between $0$ and $t$, irrespective of the value of $\alpha$, so $\hat{\beta}_n$ is equal to $0$ for all $n$ large enough (obviously, this rules out the possibility of a central limit theorem).

REMARK 9. Our estimator is based on the count of "big" increments, although we are interested in the properties of the "small" jumps of $X$, which are those governing the index $\beta$. This is because the behavior of sums of the squares of the small increments behave as described in (18), and other powers smaller than 2 are also driven by the "Wiener part" of $X$ and do not provide insight on the small jumps. Perhaps considering sums of powers bigger than 2 for "small" increments would provide an alternative means of constructing estimators of $\beta$, but we did not consider this possibility here.

2.4. *Properties of the estimators.* Our first result states that our estimators estimate $\beta$ on the (random) set $\{\bar{A}_t > 0\}$, where the jump activity index is $\beta$, and we can state the following rate of convergence.

THEOREM 1. *Let* $0 < \alpha < \alpha'$, $0 < \varpi < 1/2$ *and* $t > 0$. *Under Assumptions* 1 *and* 2, *we have* $\hat{\beta}'_n(t, \varpi, \alpha, \alpha') \xrightarrow{\mathbb{P}} \beta$ *on the set* $\{\bar{A}_t > 0\}$. *Moreover, if*

$$\chi = \chi(\beta, \gamma, \beta', \varpi)$$

(19)
$$= (\varpi\gamma) \wedge \frac{1 - \varpi\beta}{3} \wedge \frac{\varpi(\beta - \beta')}{1 + \beta'} \wedge \frac{1 - 2\varpi}{2} \wedge \frac{\varpi\beta}{2},$$

then the estimators $\hat{\beta}_n(t, \varpi, \alpha, \alpha')$ are $\Delta_n^{\chi-\varepsilon}$-rate consistent for any $\varepsilon > 0$ on the set $\{\bar{A}_t > 0\}$, in the sense that the sequence of variables $(\frac{1}{\Delta_n^{\chi-\varepsilon}}(\hat{\beta}'_n(t, \varpi, \alpha, \alpha') - \beta))_{n \geq 1}$ is bounded in probability (or, "tight") in restriction to this set.

The number $\chi$ is positive, but it may also be very small. If we want an associated distributional result, we need stronger assumptions, which essentially implies that $\chi = \varpi\beta/2$ above, and this requires that the activity indices $\beta'$ and $\beta - \gamma$ of the "nonstable-like" part of the Lévy measure be sufficiently apart from the leading activity index $\beta$, as follows.

THEOREM 2. *Let $0 < \alpha < \alpha'$ and $t > 0$. Assume Assumptions* 1 *and* 2 *with $\beta' \in [0, \beta/(2 + \beta))$ and $\gamma > \beta/2$. Then, if $\varpi < 1/(2 + \beta) \wedge 2/(5\beta)$, and in restriction to the set $\{\bar{A}_t > 0\}$, we have the following stable convergence in law to a centered normal variable independent of $X$:*

(20)
$$\frac{1}{\Delta_n^{\varpi\beta/2}}(\hat{\beta}_n(t, \varpi, \alpha, \alpha') - \beta) \xrightarrow{\mathcal{L}-(s)} \mathcal{N}\left(0, \frac{\alpha'^\beta - \alpha^\beta}{\bar{A}_t(\log(\alpha'/\alpha))^2}\right).$$

The qualifier "in restriction to the set $\{\bar{A}_t > 0\}$" is essential in this statement. Recall that, unlike the usual convergence in law, stable convergence in law makes it possible to restrict the convergence to a subset of $\Omega$ exactly as convergence in probability does. On the complement set $\{\bar{A}_t = 0\}$, anything can happen. On that set, the number $\beta$ has no meaning as a jump activity index for $X$ on $[0, t]$.

Moreover, the stable convergence in law allows for the convergence of standardized statistics.

THEOREM 3. *Under the assumptions of Theorem* 2, *the variables*

(21)
$$\frac{\log(\alpha'/\alpha)}{\sqrt{1/U(\varpi, \alpha', \Delta_n)_t - 1/U(\varpi, \alpha, \Delta_n)_t}}(\hat{\beta}_n(t, \varpi, \alpha, \alpha') - \beta)$$

*converge stably in law, in restriction to the set $\{\bar{A}_t > 0\}$, to a standard normal variable $\mathcal{N}(0, 1)$ independent of $X$.*

These results are model-free in a sense, because the drift and the volatility processes are totally unspecified apart from Assumption 1, and the Lévy measures $F_t$ are unspecified, other than the requirements specified in Assumption 2. These three theorems will be proved in Section 8 below.

The restriction on $\varpi$ given in the statement of Theorem 2 restricts admissible values of $\varpi$ in a manner that depends on $\beta$. Since $\beta$ is unknown at this point, we must select a "universal" value of $\varpi$ that is admissible for all values of $\beta$. Not surprisingly, the most stringent value of $\varpi$ is obtained in the limit where $\beta \xrightarrow{<} 2$, yielding $\varpi = 1/5$, and this is the value we suggest for empirical applications.

We note (without proof) that a similar set of properties hold for the second estimator $\hat\beta_n'$ based on the ratio of $U's$ estimated at two different frequencies $\Delta_n$ and $2\Delta_n$, with (21) replaced by the standardized statistic

$$(22) \qquad \frac{\varpi \log 2}{\sqrt{1/U(\varpi, \alpha, 2\Delta_n)_t - 1/U(\varpi, \alpha, \Delta_n)_t}}(\hat\beta_n'(t, \varpi, \alpha) - \beta).$$

**3. Stable processes.** Here, we specialize the general results in an important special case discussing, in particular, the efficiency of the estimators of $\beta$ we propose. In the special case of stable processes, the model is fully specified parametrically, and we can compare the properties of efficient parametric estimators of $\beta$ to those of the general estimators $\hat\beta_n$.

Denote, by $Y$, a symmetric stable process with index $\beta \in (0, 1/2)$. We study the two situations where $X_t = Y_t$ (the simplest of all since there is no continuous part) and $X_t = bt + \sigma W_t + Y_t$, where $\sigma > 0$, $b \in \mathbb{R}$ and $W$ is a Brownian motion. The Lévy measure depends on a scale parameter $A > 0$ and the index $\beta$. It has the form

$$(23) \qquad \begin{aligned} F(dx) &= \frac{A\beta}{2|x|^{1+\beta}}\,dx, \quad \text{hence} \\[1mm] \overline{F}(x) &:= F([-x,x]^c) = \frac{A}{x^\beta} \qquad \text{for } x > 0. \end{aligned}$$

The law of $Y_1$ has an even density $g$ and a tail function $\overline{G}(x) = \mathbb{P}(|Y_1| > x)$ satisfying, as $x \to \infty$ [see Zolotarev (1986), Theorems 2.4.2 and Corollary 2 of Theorem 2.5.1],

$$(24) \qquad g(x) = \frac{A\beta}{2|x|^{1+\beta}} + O\left(\frac{1}{x^{1+2\beta}}\right), \qquad \overline{G}(x) = \frac{A}{x^\beta} + O\left(\frac{1}{x^{2\beta}}\right).$$

In both cases $X = Y$ and $X_t = bt + \sigma W_t + Y_t$, we obviously have Assumptions 1 and 2, with $F_t = F_t' = F$ not depending on $(\omega, t)$, and with $F_t'' = F_t''' = 0$, $\Xi = \Omega \times (0, \infty)$, $\beta' = 0$, and $f_t(x) = 0$ and, finally, $A_t(\omega) = A$, which is the constant in (23). Then, we can apply the previous results, which further hold on the whole set $\Omega$ [because here $\bar{A}_t = tA > 0$ for all $(t, \omega)$]. The results are much easier to prove in this special case, and also the requirements on $\varpi$ are significantly weaker, thus allowing for faster rates of convergence (the larger $\varpi$, the faster the convergence in Theorem 1). But these improved results are no longer model-free, since the structure of jumps is completely specified in this stable model up to the unknown parameters $A$ and $\beta$.

3.1. *The case $X = Y$.* Consider first the case where $X$ has no continuous part. Then, the general results on $\hat{\beta}_n$ can be improved to yield the following.

THEOREM 4. *Assume that $X = Y$. Let $0 < \alpha < \alpha'$ and $\varpi > 0$ and $t > 0$. Then:*

(a) *If $\varpi < 1/\beta$, the estimators $\hat{\beta}_n(t, \varpi, \alpha, \alpha')$ converge in probability to $\beta$;*

(b) *If further $\varpi < 2/(3\beta)$, we have stable convergence in law, over the whole set $\Omega$, as described in Theorems 1 (with $\bar{A}_t = tA$) and 3.*

Note that in part (a) of the theorem, the closer $\beta$ is to 2, the stronger the constraint on the truncation rate $\varpi$.

These estimators are not, however, rate-efficient. To see this, one can recall from Aït-Sahalia and Jacod (2008) that the parametric model in which one observes the values $X_{i\Delta_n}$ for $i\Delta_n \leq t$ is regular, and its Fisher information for estimating $\beta$ is asymptotically of the form

$$(25) \qquad I_n \sim \frac{\log(1/\Delta_n)}{\Delta_n} C_\beta t$$

for some constant $C_\beta$. We can thus hope for estimators that, after centering by $\beta$ and normalization by $\sqrt{\log(1/\Delta_n)}/\sqrt{\Delta_n}$ are $\mathcal{N}(0, 1/C_\beta)$, and, in fact, the MLE does this.

Where is the loss of efficiency coming from? In order to compute our general estimators $\hat{\beta}_n$, we are forced by the presence of a continuous part in $X$ to discard a very sizeable portion of the data, which is the effect of truncating away the small increments of $X$. However, in this case, if somehow we knew from the start that there is no continuous part in $X$, then there would no longer be a need to do that. It is clear that better estimators of $\beta$ could then be constructed.

And if, further, the law of $Y$ has a fully-specified parametric form, as is the case here, then it would be possible to improve the estimators even more. In this example, we would simultaneously estimate $\beta$ and $A$, but the rates would be unchanged and it would be even more model-dependent. So this is the kind of estimator that we do not want to use, since we have no hope of extending such an estimator to the general semimartingale situation (or, in fact, even to more general Lévy processes than the stable ones).

3.2. *The case $X_t = bt + \sigma W_t + Y_t$.* We now study the situation where $Y$ is a stable process, but $X$ now also contains a continuous part. The distributional properties of the estimators follow directly in this special case. Indeed, for this model, $U(\varpi, \alpha)^n_t$ is essentially the same as, or close to, the number $V(\varpi, \alpha)^n_t$ of jumps of $Y$ that are bigger than $\alpha\Delta_n^\varpi$ in the interval $[0, t]$. But $V(\varpi, \alpha)^n_t$ is a Poisson random variable with parameter $Ct/\alpha^\beta \Delta_n^{\beta\varpi}$ where $C$ is a constant. Hence,

$$(26) \qquad \Delta_n^{\beta\varpi} V(\varpi, \alpha)^n_t \xrightarrow{\mathbb{P}} C/\alpha^\beta,$$

$$(27) \qquad \frac{1}{\Delta_n^{\beta\varpi/2}}(\Delta_n^{\beta\varpi} V(\varpi, \alpha)^n_t - C/\alpha^\beta) \xrightarrow{\mathcal{L}} \mathcal{N}(0, C/\alpha^\beta).$$

These properties carry over to $U(\varpi, \alpha)^n_t$, and this leads to the following improvement to Theorem 1.

THEOREM 5. *Assume that $X_t = bt + \sigma W_t + Y_t$. Let $0 < \alpha < \alpha'$ and $\varpi > 0$ and $t > 0$. Then:*

(a) *If $\varpi < 1/2$, the estimators $\hat{\beta}_n(t, \varpi, \alpha, \alpha')$ converge in probability to $\beta$;*

(b) *If $\varpi < 1/(2 + \beta)$, we have the stable convergences in law, over the whole set $\Omega$, as described in Theorems 1 (with $\bar{A}_t = tA$) and 3.*

The estimators $\hat{\beta}_n$ are again not rate-efficient, although they do come close. In fact, using the methods of Aït-Sahalia and Jacod (2008), we can show that Fisher's information for estimating $\beta$ at stage $n$ satisfies

$$(28) \qquad I_n \sim \frac{A(\log(1/\Delta_n))^{2-\beta/2}}{\sigma^\beta \Delta_n^{\beta/2}} C'_\beta t$$

for another constant $C'_\beta$. Furthermore, in the (partial) statistical model where we observe the increments provided, they are bigger than $\alpha\Delta_n^\varpi$ and discard all others (here $\alpha > 0$ and $0 < \varpi < 1/2$), Fisher's information now satisfies

$$(29) \qquad I_n \sim \frac{A(1-\varpi)^2(\log(1/\Delta_n))^2}{\alpha^\beta \Delta_n^{\varpi\beta}} C''_\beta t.$$

So, our general estimators are almost [up to a $\log(1/\Delta_n)$ factor] rate-efficient for the partial parametric statistical model. As for the "complete" model, the rate approaches the true rate by taking $\varpi$ close to $1/2$, but we cannot take $\varpi$ bigger than $1/(2 + \beta)$, and since, in practice, $\beta$ is unknown other than being less than 2, a universal choice may be $\varpi = 1/4$, which is less stringent than the choice $\varpi = 1/5$ required in the general case.

**4. General Lévy processes.** Let us now consider the case where $X$ is a general Lévy process. Its characteristics are of the form (2) with $b_t = b$, $\sigma_t = \sigma$ and $F_t = F$ deterministic and not depending on $t$. Then, Assumption 1 holds. As to Assumption 2, it may or may not hold, but if it does it takes a slightly simpler form because then everything is independent of $(\omega, t)$. In particular, $\bar{A}_t = At$ for some constant $A > 0$. The two Theorems 1 and 3 hold without modification, except that either $\{\bar{A}_t > 0\} = \Omega$ for all $t > 0$, or $\{\bar{A}_t > 0\} = \varnothing$ for all $t$, in which case those theorems are void of content.

What is important here, though, is that those results fail when the assumptions we made are not satisfied, even with such a simple probabilistic structure for $X$.

In order to see why Assumption 2 is needed, let us consider a simpler but closely related statistical model. More precisely, suppose that we observe all "big" jumps of $X$ up to time $t$; that is, $\Delta X_s$ with $|\Delta X_s| > \alpha\Delta_n^\varpi$ for all $s \le t$. A priori, this should give us more information on the Lévy measure than the original observation scheme where only increments (as opposed to jumps) are observed and only those bigger than $\alpha\Delta_n^\varpi$ are taken into consideration.

In this statistical setting, the estimators (16) have no meaning, but we can replace $\hat{\beta}_n$ with

$$(30) \qquad \overline{\beta}_n(t, \varpi, \alpha, \alpha') = \frac{\log(\overline{U}(\varpi, \alpha)_t^n / \overline{U}(\varpi, \alpha')_t^n)}{\log(\alpha'/\alpha)},$$

where we have set

$$(31) \qquad \overline{U}(\varpi, \alpha)_t^n = \sum_{s \le t} 1_{\{|\Delta X_s| > \alpha\Delta_n^\varpi\}}.$$

The estimators $\overline{\beta}_n$ are of course only virtual since there is no hope of actually observing the exact jumps of the process. But, in the rest of this section, we study the behavior of the estimators $\overline{\beta}_n$ in order to gain some insight on the necessity of making a restrictive assumption on the Lévy measure $F_t$ if one is to estimate $\beta$. We will see that such an assumption is needed even under these idealized circumstances. Set

$$(32) \qquad \gamma_n(\varpi, \alpha) = \overline{F}(\alpha\Delta_n^\varpi).$$

LEMMA 1. *Let*

$$(33) \qquad M(\varpi, \alpha)_t^n = \frac{1}{\sqrt{\gamma_n(\varpi, \alpha)}}(\overline{U}(\varpi, \alpha)_t^n - \gamma_n(\varpi, \alpha)t).$$

*Then:*

(a) *Each sequence of processes $M(\varpi, \alpha)^n$ converges stably in law to a standard Wiener process, independent of $X$;*

(b) *If $\alpha < \alpha'$, all limit points of the sequence $\gamma_n(\varpi, \alpha')/\gamma_n(\varpi, \alpha)$ are in $[0, 1]$. Further, if this sequence converges to $\gamma$, then the pair $(M(\varpi, \alpha)^n, M(\varpi, \alpha')^n)$ of processes converges stably in law to a process $(\overline{W}, \overline{W}')$, which is independent of $X$ and a 2-dimensional Wiener process with unit variances 1 and unit covariance $\sqrt{\gamma}$.*

Proof. The processes $M^n = M(\varpi, \alpha)^n$ and $M'^n = M(\varpi, \alpha')^n$ are Lévy processes and martingales, with jumps going uniformly to 0, and with predictable brackets

$$\langle M^n, M^n \rangle_t = \langle M'^n, M'^n \rangle_t = t, \qquad \langle M^n, M'^n \rangle_t = \frac{\sqrt{\gamma_n(\varpi, \alpha')}}{\sqrt{\gamma_n(\varpi, \alpha)}} t.$$

Observe, also, that $\alpha' \Delta_n^\varpi \geq \alpha \Delta_n^\varpi$; hence, $\gamma_n(\varpi, \alpha') \leq \gamma_n(\varpi, \alpha)$. The remaining results then follow [see Jacod and Shiryaev (2003), Chapter VII]. □

Theorem 6. *If* $\alpha' > \alpha$ *and* $\frac{\gamma_n(\varpi, \alpha')}{\gamma_n(\varpi, \alpha)} \to \gamma \in [0, 1]$, *then the sequence*

$$(34) \qquad \sqrt{\gamma_n(\varpi, \alpha')} \left( \overline{\beta}_n(t, \varpi, \alpha, \alpha') - \frac{\log(\gamma_n(\varpi, \alpha)/\gamma_n(\varpi, \alpha'))}{\log(\alpha'/\alpha)} \right)$$

*converges stably in law to an* $\mathcal{N}(0, \frac{1-\gamma}{t(\log(\alpha'/\alpha))^2})$ *variable, independent of* $X$.

This result is a simple consequence of the previous lemma, and its proof is the same as the proof of Theorem 1 once the CLT for the processes $U(\varpi, \alpha)^n$ is established, which we will do later.

So, the situation seems generally hopeless. These estimators are not even consistent for estimating the activity index $\beta$ of $F$ because of bias, and to remove the bias we have to know the ratio $\gamma_n(\varpi, \alpha')/\gamma_n(\varpi, \alpha)$ (or at least its asymptotic behavior in a precise way), and, further, there is no CLT if this ratio does not converge (a fact which we do not know a priori, of course).

The major difficulty comes from the possible erratic behavior of $\overline{F}$ near 0. Indeed, we have (6) with $\beta$ instead of $\beta_t^i$, but there are Lévy measures $F$ satisfying this, and such that for any $r \in (0, \beta)$ we have $x_n^r \overline{F}(x_n) \to 0$ for a sequence $x_n \to 0$ (depending on $r$, of course). If $F$ is such, the sequence $\gamma_n(\varpi, \alpha')/\gamma_n(\varpi, \alpha)$ may have the whole of $[0, 1]$ as limit points, depending on the parameter values $\varpi, \alpha, \alpha'$, and in a completely uncontrolled way for the observer.

So, we need some additional assumption on $F$. Let us consider two assumptions (the second one is stronger than the first one).

Assumption 3. $\overline{F}$ is regularly varying at 0, with index $\beta \in (0, 2)$.

Assumption 4. We have

$$(35) \qquad \overline{F}(x) = \frac{A}{x^\beta} + o\left( \frac{1}{x^{\beta/2}} \right)$$

as $x \to 0$, for some $A > 0$.

Theorem 7. (a) *Under Assumption 3, we have* $\overline{\beta}_n(t, \varpi, \alpha, \alpha') \xrightarrow{\mathbb{P}} \beta$.

(b) *Under Assumption 4, the variables* $\Delta_n^{-\varpi\beta}(\overline{\beta}_n(t,\varpi,\alpha,\alpha') - \beta)$ *converge stably in law to a* $\mathcal{N}(0, \frac{\alpha'^{\beta} - \alpha^{\beta}}{tA^2(\log(\alpha'/\alpha))^2})$ *variable, independent of* $X$.

PROOF. Assumption 3 implies that $\gamma_n(\varpi,\alpha) \to \infty$ and $\gamma_n(\varpi,\alpha)/\gamma_n(\varpi,\alpha') \to (\alpha'/\alpha)^{\beta}$, so the previous theorem yields (a). Assumption 4 clearly implies

$$\sqrt{\gamma_n(\varpi,\alpha)} \frac{\log(\gamma_n(\varpi,\alpha)/\gamma_n(\varpi',\alpha'))}{\log(\alpha'/\alpha)} \to \beta,$$

and also $\gamma_n(\varpi,\alpha) \sim A/\alpha^{\beta}\Delta_n^{\varpi\beta}$, so (b) follows again from the previous theorem. □

It may of course happen that Assumption 3 or 4 fail and nevertheless the conclusions of the previous theorem hold for a particular choice of the parameters $(\varpi,\alpha,\alpha')$ or for a particular choice of the sequence $\Delta_n$. But, in view of Theorem 6 and of the previous proof, these assumptions are *necessary* if we want those conclusions to hold *for all choices* of $(\varpi,\alpha,\alpha')$.

Now, coming back to the original realistic problem, for which only increments of $X$ are observed. Assumption 2, when $F_t(\omega,dx) = F(dx)$ for all $(\omega,t)$, is obviously stronger than Assumption 4, but not much more. The need of stronger assumptions for the original problem comes from the fact that although when we observe a "large" increment $\Delta_i^n X$ it is with a high probability almost equal to a "large" jump. Nevertheless, the observation of this jump is blurred by the Brownian component and also by a sum of very small jumps. This fact is also the reason why we need some restriction on $\varpi$ for the original problem, whereas, here, $\varpi$ can be arbitrarily large.

**5. Small sample bias correction.** By construction, we are forced by the presence of a continuous semimartingale to rely on a small fraction of the sample (i.e., those increments larger than $\alpha\Delta_n^{\varpi}$) for the purpose of estimating $\beta$. As a result, the effective sample size utilized by the estimator $\hat{\beta}_n$ is small, even if we sample at a relatively high frequency. This situation calls for an analysis of the small sample behavior of the estimator.

Such a small sample analysis is out of reach in general but it can be carried out explicitly for the model $X_t = \sigma W_t + \theta Y_t$ studied in Section 3, where $Y$ is a symmetric $\beta$-stable process and $W$ is a Wiener process. Let $g$ denote the density of $Y_1$. Here, the process $Y$ is standardized by $\mathbb{E}(e^{iuY_t}) = e^{-t|u|^{\beta}/2}$, so the limit $\beta \to 2$ corresponds to the standard normal density $\phi$.

One additional step in the expansion (24) yields, as $x \to +\infty$,

$$(36) \qquad g(x) = \frac{c_{\beta}}{x^{\beta+1}} + \frac{d_{\beta}}{x^{2\beta+1}} + O\left(\frac{1}{x^{3\beta+1}}\right)$$

and, for the tail of the distribution,

$$(37) \quad \overline{G}(x) = \mathbb{P}(|Y_1| > x) = 2 \int_x^{+\infty} g(v) \, dv = \frac{2c_\beta}{\beta x^\beta} + \frac{d_\beta}{\beta x^{2\beta}} + O\left(\frac{1}{x^{3\beta}}\right),$$

where the coefficients of the expansion are

$$(38) \quad c_\beta = \frac{\Gamma(\beta+1)}{2\pi} \sin\left(\frac{\pi\beta}{2}\right) \quad \text{and} \quad d_\beta = -\frac{\Gamma(2\beta+1)}{8\pi} \sin(\pi\beta).$$

This parametrization corresponds in terms of the general notation of the paper to

$$(39) \quad A_t = A = 2\theta^\beta c_\beta/\beta.$$

Now, consider the tail probability $P_n$ at the cutoff level $\alpha\Delta_n^\varpi$. The probability $P_n$ determines the limiting behavior of the $U$s, since $U(\varpi,\alpha)_t^n \sim (t/\Delta_n)P_n$. We have

$$\begin{aligned} P_n &= 2 \int_{-\infty}^{+\infty} \int_{\alpha\Delta_n^\varpi}^{+\infty} \frac{1}{\theta\Delta_n^{1/\beta}} g\left(\frac{x-y}{\theta\Delta_n^{1/\beta}}\right) dx \, \frac{1}{\sigma\Delta_n^{1/2}} \phi\left(\frac{y}{\sigma\Delta_n^{1/2}}\right) dy \\ &= 2 \int_{-\infty}^{+\infty} \int_{(\alpha/\theta)\Delta_n^{\varpi-1/\beta}(1-(\sigma/\alpha)\Delta_n^{1/2-\varpi}u)}^{+\infty} g(v) \, dv \, \phi(u) \, du \\ &= \int_{-\infty}^{+\infty} \overline{G}\left(\frac{\alpha}{\theta}\Delta_n^{\varpi-1/\beta}\left(1 - \frac{\sigma}{\alpha}\Delta_n^{1/2-\varpi}u\right)\right) \phi(u) \, du. \end{aligned}$$

So, with

$$\begin{aligned} (40) \quad &\left(1 - \frac{\sigma}{\alpha}\Delta_n^{1/2-\varpi}u\right)^{-\beta} \\ &= 1 + \frac{u\beta\sigma\Delta_n^{1/2-\varpi}}{\alpha} + \frac{u^2\beta(\beta+1)\sigma^2\Delta_n^{1-2\varpi}}{2\alpha^2} + O(\Delta_n^{3/2-3\varpi}) \end{aligned}$$

and $\int_{-\infty}^{+\infty} u\phi(u) \, du = 0$ and $\int_{-\infty}^{+\infty} u^2\phi(u) \, du = 1$, we see, from (37), that

$$\begin{aligned} (41) \quad P_n &= \frac{2c_\beta\theta^\beta}{\beta\alpha^\beta}\Delta_n^{1-\varpi\beta}\left(1 + \frac{\beta(\beta+1)\sigma^2}{2\alpha^2}\Delta_n^{1-2\varpi} + \frac{d_\beta\theta^\beta}{2c_\beta\alpha^\beta}\Delta_n^{1-\varpi\beta}\right) \\ &\quad + \text{smaller terms.} \end{aligned}$$

The behavior of $P_n$ suggested by the leading term in the expression $(2c_\beta\theta^\beta/(\beta\alpha^\beta))\Delta_n^{1-\varpi\beta}$ is the one we have used to define the estimator $\hat{\beta}_n$ in (16) by exploiting the dependence of that leading term on $\alpha$. The first correction term $(\beta(\beta+1)\sigma^2/(2\alpha^2))\Delta_n^{1-2\varpi}$ in (41) is due to the interaction between the Wiener and the stable processes, while the second $(d_\beta\theta^\beta/(2c_\beta\alpha^\beta))\Delta_n^{1-\varpi\beta}$ is due to the more accurate approximation of the tail of the stable process in (36) compared to the leading order term in (24).

To understand intuitively the need for the first term, suppose that the cutoff level corresponds to seven standard deviations of the continuous part of the semimartingale. There is very little probability that the Wiener process alone will generate an increment that large. On the other hand, when we count the increments due to the jump process alone, we are missing increments of the sum of the continuous and discontinuous parts where, say, the Wiener process is responsible for a one standard deviation move, and the jump process for a six standard deviation move, of the same sign. We are also missing increments where the jump process gives an eight standard deviation move and the Wiener process a one standard deviation move, of the opposite sign. The two effects partly compensate each other and indeed the term in $u$ in (40) leads to an integral whose value is zero. But the next effect, in $u^2$, leads to a net increase in the total number of increments that are larger than the cutoff when the interaction between the Wiener and jump processes is accounted for.

Asymptotically, the first of the two correcting terms in (41) is the largest, since $1 \gg \Delta_n^{1-2\varpi} \gg \Delta_n^{1-\varpi\beta}$, but in small samples a large value of the scaling parameter $\theta$ relative to $\sigma$ can make their magnitudes comparable. Using (41) at two different values $\alpha$ and $\alpha'$, we obtain

$$
\begin{aligned}
(42) \quad \hat{\beta}_n \sim \beta + \frac{1}{\log(\alpha'/\alpha)} \bigg\{ &\frac{\beta(\beta+1)\sigma^2}{2} \Big( \frac{1}{\alpha^2} - \frac{1}{\alpha'^2} \Big) \Delta_n^{1-2\varpi} \\
&+ \frac{d_\beta \theta^\beta}{2c_\beta} \Big( \frac{1}{\alpha^\beta} - \frac{1}{\alpha'^\beta} \Big) \Delta_n^{1-\varpi\beta} \bigg\}.
\end{aligned}
$$

This suggests a small sample bias correction for the estimator $\hat{\beta}_n$ obtained by subtracting an estimator of the two correction terms on the right hand side of (42) from $\hat{\beta}_n$. As we will see in simulations below, the two correction terms are quite effective in practice.

Further, we note that the two correction terms in (42), of respective orders $\Delta_n^{1-2\varpi}$ and $\Delta_n^{1-\varpi\beta}$, are asymptotically negligible at the rate $\Delta_n^{-\varpi\beta/2}$ at which the central limit occurs. This is due to the restrictions on the choice of $\varpi$ imposed by Theorem 2. Consequently, the bias-corrected estimator has the same asymptotic distribution as the original estimator.

More generally, we have

$$
\begin{aligned}
(43) \quad \hat{\beta}_n \sim \beta + \frac{1}{\log(\alpha'/\alpha)} \bigg\{ &\frac{\int_0^t A_s \sigma_s^2 \, ds}{\bar{A}_t} \frac{\beta(\beta+1)}{2} \Big( \frac{1}{\alpha^2} - \frac{1}{\alpha'^2} \Big) \Delta_n^{1-2\varpi} \\
&+ \frac{\int_0^t A_s^2 \, ds}{\bar{A}_t} \frac{\beta d_\beta}{4c_\beta^2} \Big( \frac{1}{\alpha^\beta} - \frac{1}{\alpha'^\beta} \Big) \Delta_n^{1-\varpi\beta} \bigg\}.
\end{aligned}
$$

To implement the bias correction in practice, we need to estimate the terms $(1/\bar{A}_t)\int_0^t A_s \sigma_s^2 \, ds$ and $(1/\bar{A}_t)\int_0^t A_s^2 \, ds$. In the case of a stable symmetric process, $A_s = A$ and so $(1/\bar{A}_t)\int_0^t A_s^2 \, ds = A = 2\theta^\beta c_\beta/\beta$. We can then replace $(1/\bar{A}_t)\int_0^t A_s \sigma_s^2 \, ds$ by $(1/t)\int_0^t \sigma_s^2 \, ds$ and use any standard estimator of the integrated volatility. In general, we have

$$
\begin{aligned}
(44) \quad U(\varpi,\alpha)_t^n \sim \frac{1}{\alpha^\beta}\Delta_n^{-\varpi\beta}\Big(&\bar{A} + \frac{1}{\alpha^2}\overline{A\sigma^2}G_1(\beta)\Delta_n^{1-2\varpi} \\
&+ \frac{1}{\alpha^\beta}\overline{A^2}G_2(\beta)\Delta_n^{1-\varpi\beta}\Big),
\end{aligned}
$$

$$
(45) \qquad \Delta_n^{\varpi\beta}U(\varpi,\alpha)_t^n \sim a_0\frac{1}{\alpha^\beta} + a_1\frac{1}{\alpha^{2+\beta}} + a_2\frac{1}{\alpha^{2\beta}},
$$

where $\bar{A} = \int_0^t A_s \sigma_s^2 \, ds$, $\overline{A\sigma^2} = \int_0^t A_s \sigma_s^2 \, ds$ and $\overline{A^2} = \int_0^t A_s^2 \, ds$. We can estimate the unknown coefficients $a_0$, $a_1$ and $a_2$ in expression (45) by a straightforward linear regression of $\Delta_n^{\varpi\beta}U(\varpi,\alpha)_t^n$ on $1/\alpha^\beta$, $1/\alpha^{2+\beta}$ and $1/\alpha^{2\beta}$. For the purpose of running that regression, we use different cutoff levels $\alpha$ and compute the corresponding number of increments exceeding that level, $U(\varpi,\alpha)_t^n$ and the first-stage estimate of $\beta$. Given estimates of the regression coefficients, we have a generalized bias correction procedure based on subtracting, from $\hat{\beta}_n$, the terms on the right-hand side of

$$
(46) \qquad \hat{\beta}_n - \beta \sim \frac{1}{\log(\alpha'/\alpha)}\left\{\frac{a_1}{a_0}\left(\frac{1}{\alpha^2}-\frac{1}{\alpha'^2}\right) + \frac{a_2}{a_0}\left(\frac{1}{\alpha^\beta}-\frac{1}{\alpha'^\beta}\right)\right\}
$$

evaluated at the regression estimates of $a_0$, $a_1$ and $a_2$.

**6. Monte Carlo simulations.** We now report simulation results documenting the finite sample performance of the estimator $\hat{\beta}_n$ in finite samples. We calibrate the values to be realistic for a very liquid stock. We use an observation length of $T = 1$ day, consisting of 6.5 hours of trading (i.e., $n = 23,400$ seconds).

The averages and standard deviations of the estimator $\hat{\beta}_n$, which is based on two different levels of truncation $\alpha$ and $\alpha'$, are reported in Table 1 for various values of $\beta$ up to 1.5 and include a continuous (Brownian) part. The table reports the results of 5000 simulations. The data generating process is the stochastic volatility model $dX_t = \sigma_t \, dW_t + \theta \, dY_t$, with $\sigma_t = v_t^{1/2}$, $dv_t = \kappa(\eta - v_t)\,dt + \gamma v_t^{1/2}\,dB_t + dJ_t$, $E[dW_t\,dB_t] = \rho\,dt$, $\eta^{1/2} = 0.25$, $\gamma = 0.5$, $\kappa = 5$, $\rho = -0.5$. $J$ is a compound Poisson jump process with jumps that are uniformly distributed on $[-30\%, 30\%]$ and $X_0 = 1$. The jump process $Y$ is either a $\beta$-stable process with $\beta = 1.5$, 1.25, 1.0, 0.75, 0.5 and 0.25, or a compound Poisson process (which has finite activity and is marked $\beta = 0$ in the table) with fixed jump size 0.10. The estimator is implemented with $\alpha = 5\eta$, $\alpha' = 10\eta$ and $\varpi = 0.20$. Given $\eta$ and $\alpha$, the scale parameter $\theta$ (or equivalently $A$) of the stable process in simulations is calibrated to deliver the various values of the tail probability $\mathbb{P}(|\Delta Y_t| \geq \alpha \Delta_n^{\varpi})$ reported in the columns of the table; for the Poisson process, it is the value of the arrival rate parameter $\lambda$ that is set to generate the desired level of jump tail probability.

TABLE 1
*Monte Carlo simulations of the estimator $\hat{\beta}_n$ based on two levels of truncation for*
*$\beta$-stable processes and a compound Poisson process ($\beta = 0$)*

| Sampling $\Delta_n$ tail probability | | 1 sec 0.25% | 1 sec 0.5% | 1 sec 1.0% | 1 sec 2.5% | 5 sec 1.0% |
|---|---|---|---|---|---|---|
| $\beta = 1.5$ | Sample mean | 1.52 | 1.51 | 1.50 | 1.52 | 1.53 |
| | Sample stdev | (0.26) | (0.18) | (0.13) | (0.08) | (0.25) |
| | Asymp stdev | (0.26) | (0.18) | (0.13) | (0.08) | (0.24) |
| $\beta = 1.25$ | Sample mean | 1.27 | 1.26 | 1.25 | 1.26 | 1.27 |
| | Sample stdev | (0.23) | (0.16) | (0.11) | (0.07) | (0.19) |
| | Asymp stdev | (0.23) | (0.16) | (0.11) | (0.07) | (0.19) |
| $\beta = 1.0$ | Sample mean | 1.01 | 1.01 | 1.00 | 1.00 | 1.01 |
| | Sample stdev | (0.19) | (0.14) | (0.10) | (0.06) | (0.14) |
| | Asymp stdev | (0.19) | (0.14) | (0.10) | (0.06) | (0.14) |
| $\beta = 0.75$ | Sample mean | 0.76 | 0.76 | 0.75 | 0.75 | 0.76 |
| | Sample stdev | (0.16) | (0.11) | (0.08) | (0.05) | (0.11) |
| | Asymp stdev | (0.16) | (0.11) | (0.08) | (0.05) | (0.11) |
| $\beta = 0.5$ | Sample mean | 0.51 | 0.50 | 0.50 | 0.50 | 0.50 |
| | Sample stdev | (0.13) | (0.09) | (0.06) | (0.04) | (0.08) |
| | Asymp stdev | (0.13) | (0.09) | (0.06) | (0.04) | (0.08) |
| $\beta = 0.25$ | Sample mean | 0.25 | 0.25 | 0.25 | 0.25 | 0.25 |
| | Sample stdev | (0.09) | (0.06) | (0.04) | (0.04) | (0.05) |
| | Asymp stdev | (0.09) | (0.06) | (0.04) | (0.04) | (0.05) |
| $\beta = 0$ | Sample mean | 0.01 | 0.01 | 0.01 | 0.01 | 0.02 |
| | Sample stdev | (0.02) | (0.01) | (0.007) | (0.005) | (0.01) |

In each row, the top number is the average value of the estimator $\hat{\beta}_n$ across the simulations, after inclusion of the bias correction discussed in Section 5, while the number below, in parentheses, is the standard deviation of the estimator across the same simulations. The third number in parentheses is the estimated asymptotic standard error based on the limiting distribution given in the sections above. A higher tail probability in the columns has the effect of generating more increments from the jump process that exceed the cutoff level, which makes more observations available and correspondingly reduces the standard deviation of the estimates.

As the results show, $\hat{\beta}_n$ picks up on average fairly accurately the true value of $\beta$. As $\beta$ gets too close to 2, the $\beta$-stable jump process starts to approximate

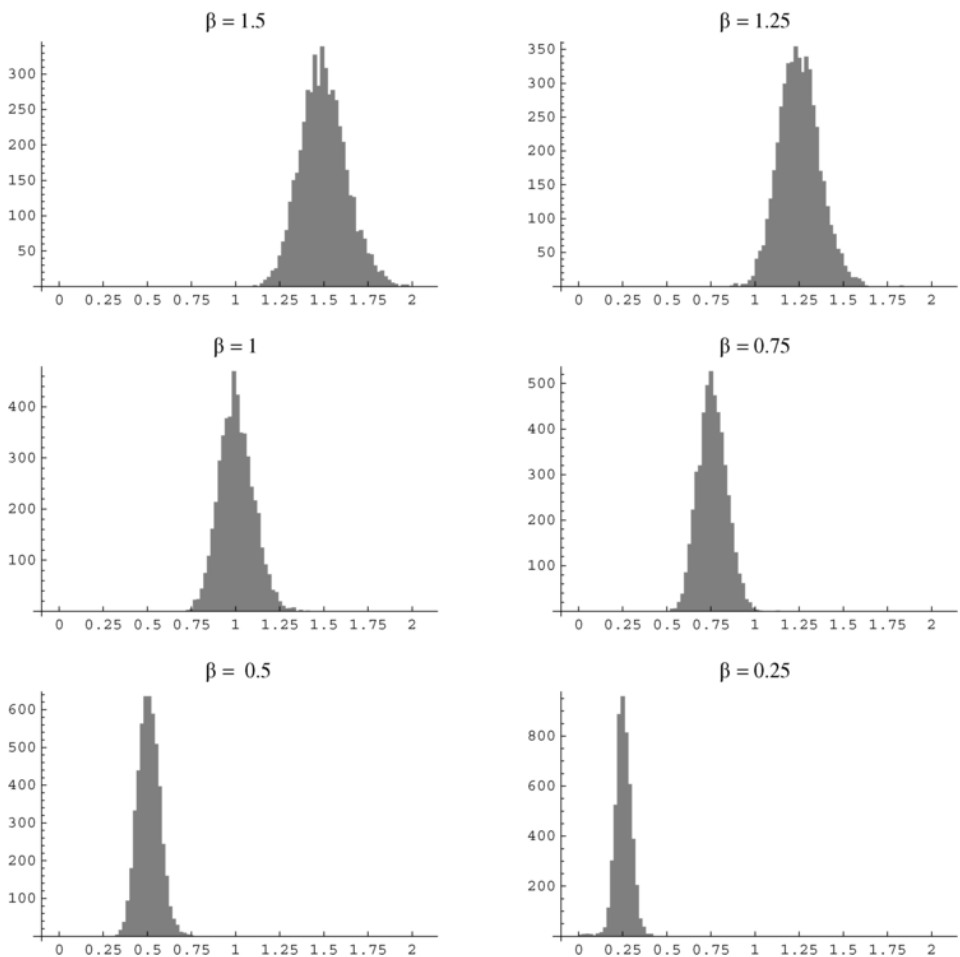

Fig. 1. *Monte Carlo distributions of the estimator $\hat{\hat{\beta}}_n$ based on two levels of truncation for $\beta$-stable processes $(0 < \beta < 2)$.*

too closely the behavior of the Brownian motion and the performance of the estimator deteriorates.

Further simulations (not reported to save space) suggest that the estimator is not overly sensitive to the selection of the truncation levels $(\alpha, \alpha')$ within a reasonable range. Histograms of the distribution of the estimator $\hat{\beta}_n$ are shown in Figure 1 for the same values of $\beta$; the figures are based on the 1% level of jump tail probability. The histogram reports the raw, unstandardized, values of $\hat{\beta}_n$, which are not expected to be asymptotically normal unlike the standardized versions.

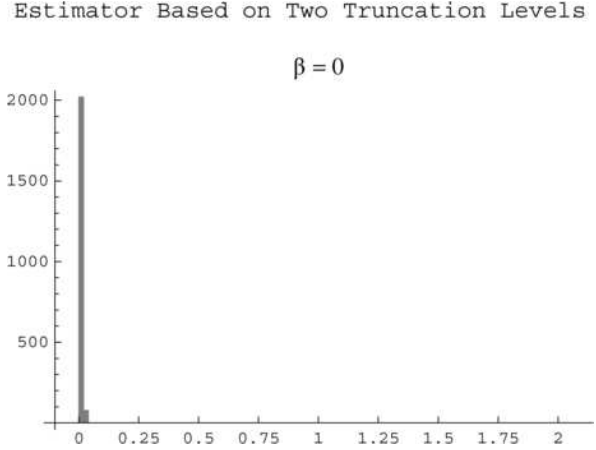

Fig. 2. *Monte Carlo distribution of the estimator $\hat{\beta}_n$ based on two levels of truncation for a compound Poisson process ($\beta = 0$).*

A special mention should be made about the compound Poisson case, which corresponds to $\beta = 0$ and does not satisfy Assumption 2. As discussed in Remark 8, there is no central limit theorem in this case, and $\hat{\beta}_n$ should be equal to 0 for $n$ large enough. Intuitively, when there is a small number or large jumps, the same large increments remain in the sample at the two truncation levels $\alpha$ and $\alpha'$, and the ratio of $U$ evaluated at $\alpha$ to $U$ evaluated at $\alpha'$ in (16) is equal to 1. This is what happens in simulations in the vast majority of cases, as shown by the histogram for $\beta = 0$ reported in Figure 2.

The asymptotic distribution is an accurate guide for the small samples as shown in the standardized distributions in Figure 3. The estimator in the histograms is standardized according to the asymptotic distribution given in Theorem 3, and the solid curve in the figures is the limiting $\mathcal{N}(0, 1)$ density. As the figures show, the asymptotic distribution is a fairly accurate guide for the small samples. This is in spite of the relatively small number of (large) increments that are effectively used by the estimator, combined with the facts that some large increments are kept, even though they may not have contained a large jump, or, conversely, smaller increments may have contained two or more large, cancelling, jumps, or the Wiener process may have combined with the pure jump process to produce a larger increment. Asymptotically, these effects do not show up at the leading order in $\Delta_n$ but are present in small samples and appear to be effectively captured by the bias correcting term.

Finally, we compare in simulations the performance of the two estimators $\hat{\beta}_n$ (based on two truncation levels) and $\hat{\beta}'_n$ (based on two sampling frequencies) with the same experiment design as above. The sampling frequency is $\Delta_n = 1$ second, the length of observation $T = 1$ day or 23,400 seconds. The tail probability is 1.0%, the middle value in Table 1. The results are shown in Table 2. The estimator $\hat{\beta}'_n$ tends to have a larger standard deviation than $\hat{\beta}_n$ and is slightly biased. For these reasons, we have focused on the estimator $\hat{\beta}_n$ and emphasize its use in the empirical application that follows.

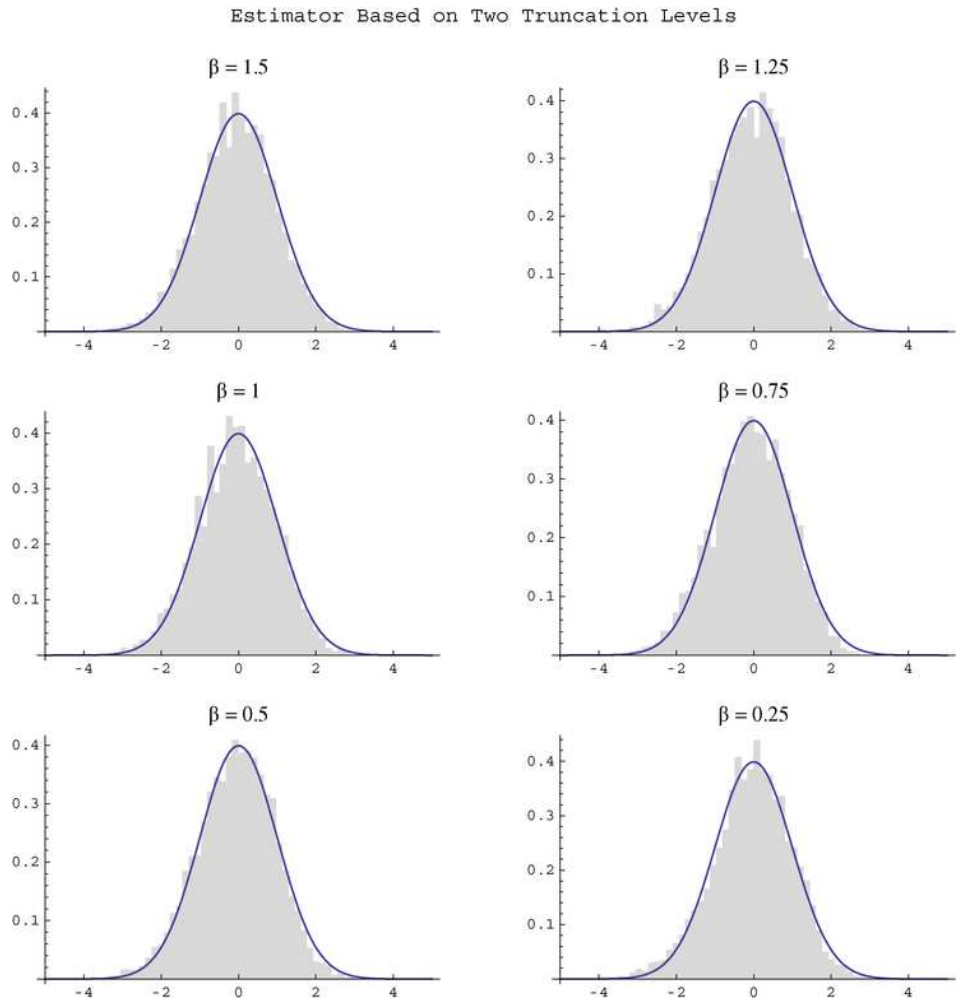

FIG. 3.   *Standardized Monte Carlo and asymptotic distributions of the estimator $\hat{\beta}_n$ based on two levels of truncation for $\beta$-stable processes $(0 < \beta < 2)$.*

**7. Empirical application.** We now implement the estimator $\hat{\beta}_n$ for the two most actively traded stocks in the Dow Jones Industrial Average index, Intel (INTC) and Microsoft (MSFT), and each trading day in 2006. The data source is the TAQ database. Each day, we collect all transactions on the NYSE or NASDAQ, from 9:30 am until 4:00 pm, for each one of these stocks. We sample in calendar time every 5 and 15 seconds. We use filters to eliminate clear data errors (price set to zero, etc.) and all transactions in the original record that are later corrected, cancelled or otherwise invalidated, as is standard in the empirical high frequency literature.

TABLE 2
*Comparison of the two estimators of $\beta$ in Monte Carlo simulations for $\beta$-stable processes and a compound Poisson process ($\beta = 0$)*

|  |  | Two truncation levels $\hat{\beta}_n$ | Two sampling frequencies $\hat{\beta}'_n$ |
|---|---|---|---|
| $\beta = 1.0$ | Sample mean | 1.00 | 1.00 |
|  | Sample stdev | (0.10) | (0.13) |
| $\beta = 0.5$ | Sample mean | 0.50 | 0.51 |
|  | Sample stdev | (0.06) | (0.09) |
| $\beta = 0$ | Sample mean | 0.01 | 0.03 |
|  | Sample stdev | (0.007) | (0.02) |

The two time series are plotted in Figure 4. Figure 5 contains a histogram of the tails of the unconditional densities of the log-returns from the two stocks. Comparing the two figures, we see that it is quite possible to have a standard time series plot display little evidence of large moves (Figure 5), while the tails of the distribution look substantially fatter than normal (Figure 5) as confirmed by the descriptive statistics in Table 3 for the two log-returns series. All together, this evidence points in the direction of many small, active jumps of the type that we seek to uncover using our estimator of $\beta$.

More formally, we compute the statistic $\hat{S}_n$ of Aït-Sahalia and Jacod (2009) to test for the presence of jumps in the data. Over the different sampling frequencies considered (ranging from 2 seconds to 1 minute), the largest value of the statistic $\hat{S}_n$ we obtain for the different quarters and the two stocks is 1.19. Since the asymptotic value of $\hat{S}_n$ is 1 (resp. 2) when jumps

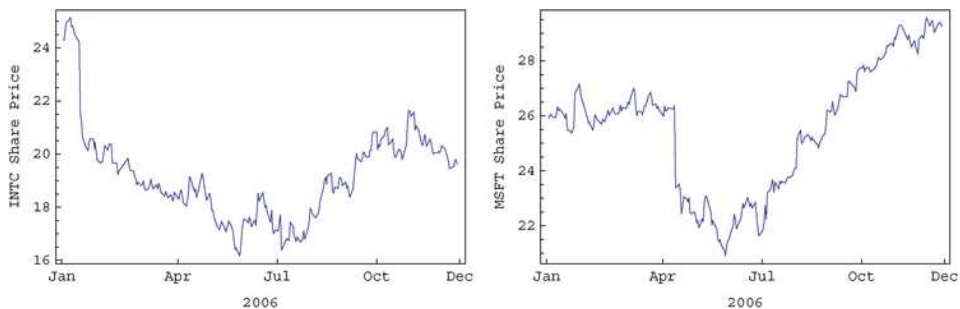

FIG. 4. *Time series of INTC and MSFT stock prices, all trading days in 2006.*

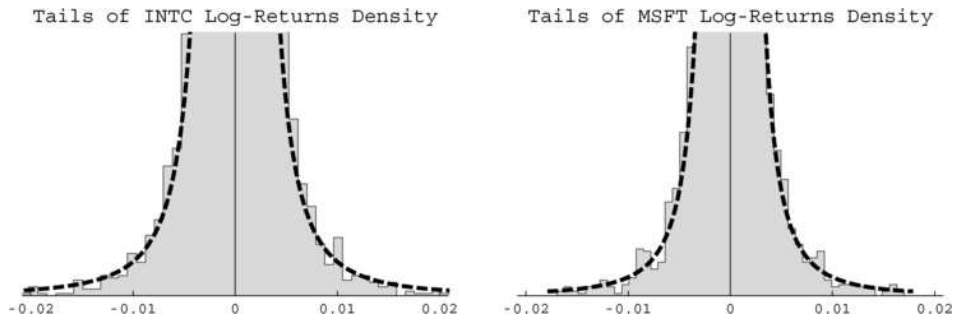

Fig. 5.   *Tails of the marginal density of INTC and MSFT log-returns, all trading days in 2006 sampled at 5 second intervals. The dashed curve is the leading part of the Lévy measure, that is $F'_t(dx)$ given in (8), with jump activity index $\beta$, estimated for each stock by $\hat{\beta}_n$.*

are present (resp. absent), this provides further evidence for the presence of jumps.

Table 4 reports the estimates produced by $\hat{\beta}_n$ for the two stocks, for the full year. For each stock and sampling intervals $\Delta_n$ of 5 and 15 seconds, respectively, the left columns in the table use a level of truncation $\alpha$, which is set to 7 estimated standard deviations of the continuous part of the process, a second level $\alpha' = 2\alpha$ and a truncation rate $\varpi = 0.20$. In light of the different possible choices of $\alpha$ and $\alpha'$, a simple variance-reducing procedure consists in averaging the estimators obtained over the different $(\alpha, \alpha')$. The columns marked AVG report the results of averaging the estimates over values of $\alpha$ ranging from 7 to 9 estimated standard deviations of the continuous part of the process and over values of $\alpha'$ ranging from 1.5 to 4 times $\alpha$, as discussed above in Remark 7. $\hat{\beta}_n$ is the estimator defined in (16). The estimator $\tilde{\beta}_n$ in the table denotes the estimator obtained by applying to $\hat{\beta}_n$ the regression-based bias-correction procedure described in Section 5 implemented using the same range of truncation levels as for the average AVG. Standard errors of the estimators are in parentheses.

TABLE 3
*Descriptive statistics for Intel and Microsoft transactions in 2006*

|          | Mean              | Stdev   | Skew | Kurt | Min    | Max   |
|----------|-------------------|---------|------|------|--------|-------|
|          |                   |         | INTC |      |        |       |
| Qtr 1    | $-2 \times 10^{-6}$ | 0.00069 | $-33$ | 5195 | $-0.109$ | 0.019 |
| Qtr 2    | $-3 \times 10^{-7}$ | 0.00070 | $-0.1$ | 136  | $-0.027$ | 0.027 |
| Qtr 3    | $6 \times 10^{-7}$  | 0.00071 | $-0.4$ | 93   | $-0.027$ | 0.022 |
| Qtr 4    | $-1 \times 10^{-7}$ | 0.00065 | $-0.5$ | 110  | $-0.021$ | 0.030 |
| All year | $-5 \times 10^{-7}$ | 0.00069 | $-8.5$ | 1430 | $-0.109$ | 0.030 |
|          |                   |         | MSFT |      |        |       |
| Qtr 1    | $3 \times 10^{-7}$  | 0.00050 | 0.1  | 190  | $-0.025$ | 0.028 |
| Qtr 2    | $-1 \times 10^{-6}$ | 0.00065 | $-49$ | 8927 | $-0.119$ | 0.021 |
| Qtr 3    | $1 \times 10^{-6}$  | 0.00053 | 7.3  | 749  | $-0.017$ | 0.051 |
| Qtr 4    | $7 \times 10^{-7}$  | 0.00051 | $-0.2$ | 388  | $-0.031$ | 0.031 |
| All year | $2 \times 10^{-7}$  | 0.00055 | $-19$ | 4715 | $-0.119$ | 0.051 |

Table 4
*Estimates of β from all Intel and Microsoft transactions in 2006*

| $\Delta_n$ | INTC | | | | MSFT | | | |
|---|---|---|---|---|---|---|---|---|
| | 5 sec | | 15 sec | | 5 sec | | 15 sec | |
| $\alpha$ | **7** | **AVG** | **7** | **AVG** | **7** | **AVG** | **7** | **AVG** |
| $\alpha'/\alpha$ | **2** | **AVG** | **2** | **AVG** | **2** | **AVG** | **2** | **AVG** |
| $\hat{\beta}_n$ | 1.43 | 1.56 | 1.76 | 1.72 | 1.69 | 1.62 | 1.60 | 1.61 |
| $\tilde{\beta}_n$ | | 1.52 | | 1.69 | | 1.60 | | 1.59 |
| | (0.04) | (0.003) | (0.05) | (0.006) | (0.05) | (0.004) | (0.05) | (0.005) |

We find evidence of infinitely active jumps. The estimated values of $\beta$ are about 1.5 for INTC and 1.6 for MSFT. These correspond also to the average values $\hat{\beta}_n$ that we obtain if we average the four quarterly estimates, and also average them over different values of $\alpha$ and $\alpha'/\alpha$. Clearly, the nature of these infinitely active jumps cannot be assessed by mere visual inspection of the time series in Figure 4, as would have been the case if only large, infrequent jumps, were present.

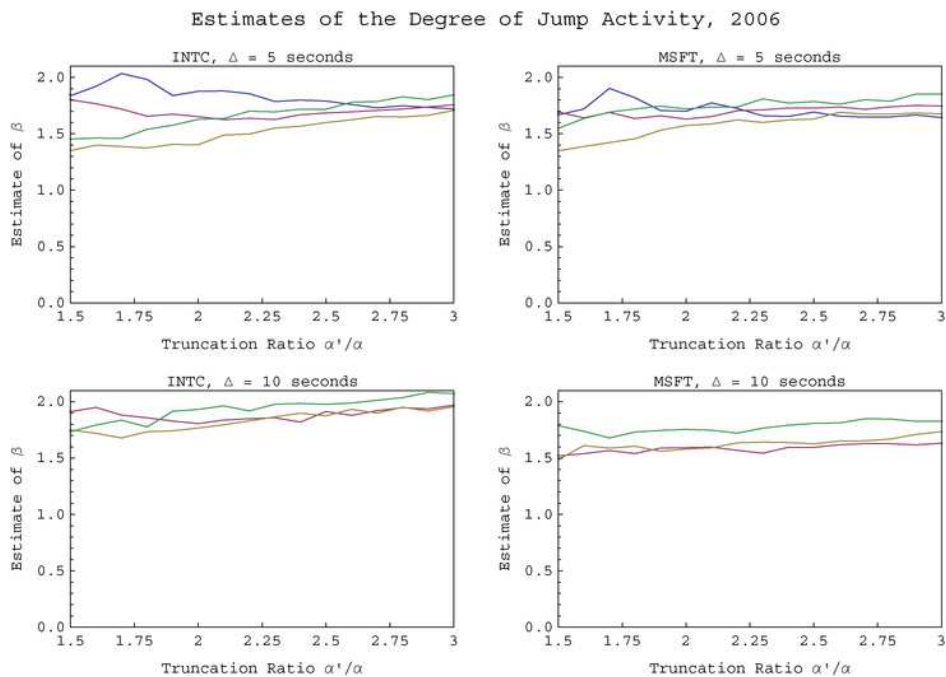

Fig. 6. *Estimates of β obtained from $\hat{\beta}_n$ at 5 and 15 seconds for Intel and Microsoft, all 2006 transactions.*

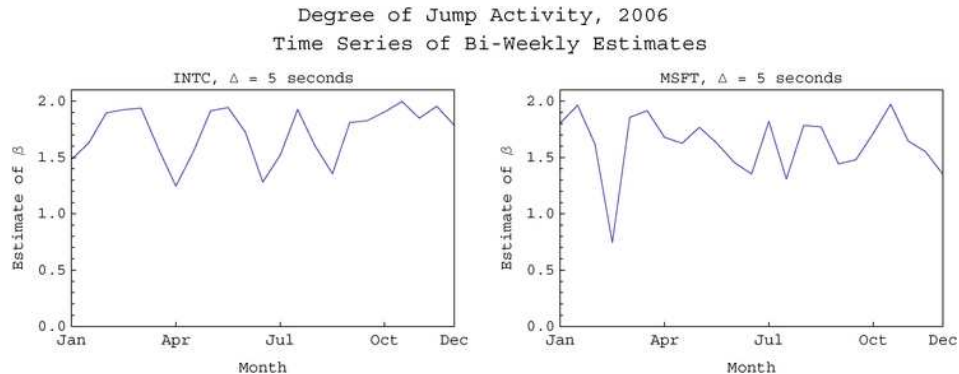

FIG. 7.  *Time series of the estimated β for Intel and Microsoft, all 2006 transactions, computed at the 5 second frequency on a bi-weekly basis.*

Figure 6 shows the values of the estimates computed over the full year for various values of $\alpha$ (the four curves on each plot correspond, resp., to 6, 7, 8 and 9 estimated standard deviations of the continuous part of the process) and values of $\alpha'$ ranging from 1.5 to 3 times $\alpha$ (the horizontal axis represents the ratio $\alpha'/\alpha$). The continuous standard deviations are computed from the small increments of the process and updated each quarter. The figure shows that, for each stock, we obtain comparable values for $\beta$ when varying $\alpha$ and $\alpha'$. The top panel of the plot shows the estimates at $\Delta_n = 5$ seconds, the bottom one at $\Delta_n = 15$ seconds. As expected, the lower frequency estimates are more variable, since they rely on a smaller number of large increments in total, and the variance of $\hat{\beta}_n$ is proportional to the difference between one over the number of increments retained at the two cutoff levels.

Figure 7 plots the time series of the estimated $\beta$ for Intel and Microsoft computed using a time interval $T$ of two weeks during year 2006 at the 5 second frequency. The estimates are the AVG estimates (averages over the values obtained for a range of $\alpha$ and $\alpha'$) discussed in Table 4. The figure shows that the bi-weekly estimates, while variable over the course of the year, tend to be clustered around the same value of 1.5 as the full-year estimates. Shorter lengths of observation than two weeks tend to produce substantially more variable estimates, due to the small number of large increments available.

Allow us one last comment regarding the high sampling frequency and its interaction with market microstructure noise in the data. There is no doubt that market microstructure noise is a concern at sampling intervals of a few seconds. However, in the present context, we are only using the large increments of the process, namely those greater than $\alpha\Delta^\varpi$. Such increments are more likely to represent a true price movement rather than noise, especially after the standard data cleaning consisting of removing price errors have been applied, whereas the noise most likely mainly corrupts the observations by a relatively small amount. Similarly, the cutoff points are beyond the level of bid/ask bounces, or the tick size of 1 cent.

Furthermore, all we are doing is counting those increments, not exploiting their magnitude. So, provided that the noise does not substantially affect the likelihood that an increment is above or below the cutoff, these estimates are unlikely to be seriously affected by the noise. But this paper represents only a first attempt at measuring the degree of jump activity. We leave to future work the development of estimators that are fully robust to market microstructure noise.

**8. Technical results and proofs.** Our basic semimartingale $X$ satisfies (2) and Assumptions 1 and 2. Note that all three Theorems 1, 2 and 3 are "local" in time, so by a standard localization procedure [see, e.g., Aït-Sahalia and Jacod (2009)] we can indeed assume a strengthened version of our assumptions, which follow.

ASSUMPTION 5. The processes $b$ and $\sigma$ are bounded by some constant $L$.

ASSUMPTION 6. We have Assumption 2 with $L_t(\omega) = L$ being constant.

In all the sequels, the letter $K$ denotes a constant that changes from line to line and may depend on $X$ and its characteristics and also on the parameters $\varpi, \alpha, \alpha'$, and we write $K_p$ if we want to emphasize its dependency on some other parameter $p$. We use the shorthand notation $\mathbb{E}_{i-1}^n$ and $\mathbb{P}_{i-1}^n$, respectively, $\mathbb{E}_t$ and $\mathbb{P}_t$ or the conditional expectation and probability with respect to $\mathcal{F}_{(i-1)\Delta_n}$, respectively $\mathcal{F}_t$.

Before proceeding, we mention a number of elementary consequences of Assumption 6, to be used many times. For all $u \in (0, 1]$ and $v \in (0, 2]$ and

$x, y \in (0, 1]$, we have [recall the notation (4)]

$$(47) \quad \begin{cases} \overline{F}_t''(x) \leq \dfrac{K}{x^{\beta'}}, \qquad \left| \overline{F}_t(x) - \dfrac{A_t}{x^{\beta}} \right| \leq \dfrac{K}{x^{(\beta-\gamma)\vee\beta'}}, \qquad \overline{F}_t(x) \leq \dfrac{K}{x^{\beta}}, \\[2mm] \displaystyle\int_{\{|x|\leq u\}} x^2 F_t(dx) \leq K u^{2-\beta}, \\[2mm] \displaystyle\int_{\{|x|>u\}} (|x|^v \wedge 1) F_t(dx) \leq \begin{cases} K_v, & \text{if } v > \beta, \\ K_v \log(1/u), & \text{if } v = \beta, \\ K_v u^{v-\beta}, & \text{if } v < \beta, \end{cases} \\[2mm] \overline{F}_t(x) - \overline{F}_t(x+y) \leq \dfrac{K}{x^{\beta}}\left(1 \wedge \dfrac{y}{x} + x^{\gamma \wedge (\beta-\beta')}\right). \end{cases}$$

8.1. *Estimates for stable processes.* In this subsection, we consider a symmetric stable process $Y$ with Lévy measure (23). For each $\delta \in (0, 1]$, we set

$$(48) \qquad Y(\delta)_t' = Y_t - \sum_{s \leq t} \Delta Y_s 1_{\{|\Delta Y_s| > \delta\}}.$$

LEMMA 2. *There is a constant $K$, depending on $(A, \beta)$, such that, for all $s > 0$,*

$$(49) \qquad \mathbb{P}(|Y(\delta)_s'| > \delta/2) \leq K s^{4/3}/\delta^{4\beta/3}.$$

PROOF. We use the notation (23) and (24). We set $\theta = s\overline{F}(\delta/2) = sA(2/\delta)^{\beta}$ and consider the processes $Y' = Y(\delta)'$, $Y'' = Y(\delta/2)'$ and $Z_t = \sum_{r \leq t} 1_{\{|\Delta Y_r| > \delta/2\}}$. Introduce, also, the sets

$$D = \{|Y_s| > \delta/2\}, \qquad D' = \{|Y_s'| > \delta/2\},$$
$$B = \{Z_s = 1\}, \qquad B' = \{Z_s = 0\}.$$

It is of course enough to prove the result for $s/\delta^{\beta}$ small, so, below, we assume $\theta \leq 1/2$.

By scaling, $\mathbb{P}(D) = \overline{G}(\delta s^{-1/\beta}/2)$, so (24) yields

$$(50) \qquad |\mathbb{P}(D) - \theta| \leq K\theta^2.$$

Next, $Z_s$ is a Poison variable with parameter $\theta \leq 1/2$. Hence,

$$(51) \qquad |\mathbb{P}(B) - \theta| \leq K\theta^2.$$

Since $Y''$ is a purely discontinuous Lévy process without drift and whose Lévy measure is the restriction of $F$ to $[-\delta', \delta']$, we have

$$(52) \qquad \mathbb{E}((Y_s'')^2) = s \int_{\{|x| \leq \delta/2\}} x^2 F(dx) \leq K\theta\delta^2.$$

The two processes $Y''$ and $Z$ are independent, and conditionally on $B$ the law of the variable $Z_s$ is $F$ restricted to $\{|x| > \delta/2\}$ and normalized by $\theta$. Hence,

$$
\begin{aligned}
\mathbb{P}(B \cap D^c) &= e^{-\theta} s \int_{\{|x| > \frac{\delta}{2}\}} F(dx) \mathbb{P}\left(|Y''_s + x| \le \frac{\delta}{2}\right) \\
&\le e^{-\theta} s \int_{\{|x| > \delta/2\}} F(dx) \mathbb{P}\left(|Y''_s| \ge |x| - \frac{\delta}{2}\right) \\
&\le s\left(F\left(\left\{\frac{\delta}{2} < |x| \le \frac{\delta}{2}(1 + \theta^{1/3})\right\}\right)\right. \\
&\qquad \left. + F\left(\left\{|x| > \frac{\delta}{2}\right\}\right) \mathbb{P}\left(|Y''_s| > \frac{\delta}{2} \theta^{1/3}\right)\right) \\
&\le K\theta\left(\theta^{1/3} + \frac{1}{\delta^2 \theta^{2/3}} \mathbb{E}((Y''_s)^2)\right) \le K\theta^{4/3},
\end{aligned}
$$
(53)

where we have used (23) and (52) for the last two inequalities.

Now, we have

$$
\mathbb{P}(D \cap B^c) = \mathbb{P}(D) - \mathbb{P}(B) + \mathbb{P}(B \cap D^c).
$$

Observe, also, that $D \cap B' = D' \cap B'$, and $D'$ and $B'$ are independent. Hence,

$$
\mathbb{P}(D') = \frac{\mathbb{P}(D' \cap B')}{\mathbb{P}(B')} = \frac{\mathbb{P}(D \cap B')}{\mathbb{P}(B')} \le \frac{\mathbb{P}(D \cap B^c)}{\mathbb{P}(B')} \le K\mathbb{P}(D \cap B^c)
$$

because $\mathbb{P}(B') = e^{-\theta} \ge e^{-1/2}$. The last two displays, plus (50) and (51) and (53) give us $\mathbb{P}(D') \le K\theta^{4/3}$. Hence, (49).  $\square$

### 8.2. *Estimates for semimartingales.* 

Below, we assume Assumptions 5 and 6 without special mention. The semimartingale $X$ can be written as

$$
X = X_0 + B + X^c + (x1_{\{|x| \le 1\}}) \star (\mu - \nu) + (x1_{\{|x| > 1\}}) \star \mu,
$$

where $X^c$ is the continuous martingale part and $\mu$ is the jump measure, $B$ and $\nu$ are as in (2), and the "$\star$" stands for the stochastic integral with respect to random measures [see Jacod and Shiryaev (2003)]. For any $\delta \in (0, 1]$ we set

$$
(54) \quad \begin{cases} X(\delta)''_t = \sum_{s \le t} \Delta X_s 1_{\{|\Delta X_s| > \delta\}}, \\ X(\delta)' = X - X(\delta)'' = X_0 + B + X^c + (x1_{\{|x| \le \delta\}}) \star (\mu - \nu) - B(\delta), \end{cases}
$$

where

$$
(55) \qquad B(\delta)_t = \int_0^t b(\delta)_s \, ds, \qquad b(\delta)_s = \int_{\{\delta < |x| \le 1\}} x F_s(dx).
$$

By Assumption 6, $F'_t$ and $F''_t$ are mutually singular. Hence, there exists a predictable subset $\Phi$ of $\Omega \times (0, \infty) \times \mathbb{R}$ such that

$$(56) \qquad \begin{cases} F'_t(\omega, \cdot) \text{ is supported by the set } \{x : (\omega, t, x) \notin \Phi\}, \\ F''_t(\omega, \cdot) \text{ is supported by the set } \{x : (\omega, t, x) \in \Phi\}. \end{cases}$$

Observe, also, that if $\beta' \leq 1$, we can set

$$(57) \qquad B(\delta)'_t = \int_0^t b(\delta)'_s \, ds, \qquad b(\delta)'_s = \int_{\{\delta < |x| \leq 1\}} x F''_s(dx)$$

and $|b(\delta)'_s| \leq K$. In this case, $B(\delta)' = (x 1_{\{|x| \leq \delta\}} 1_\Phi) \star \nu$ is of finite variation, and so is $(x 1_{\{|x| \leq \delta\}} 1_\Phi) \star \mu$. Therefore, we have the decomposition

$$(58) \qquad \begin{cases} X(\delta)' = X_0 + \widehat{X} + \overline{X}(\delta)^a + \overline{X}(\delta)^b - B(\delta), \qquad \text{where} \\ \overline{X}(\delta)^a = \begin{cases} (x 1_{\{|x| \leq \delta\}} 1_\Phi) \star \mu, & \text{if } \beta' \leq 1, \\ (x 1_{\{|x| \leq \delta\}} 1_\Phi) \star (\mu - \nu), & \text{if } \beta' > 1, \end{cases} \\ \overline{X}(\delta)^b = (x 1_{\{|x| \leq \delta\}} 1_{\Phi^c}) \star (\mu - \nu), \\ \widehat{X} = \begin{cases} B + X^c - B(\delta)', & \text{if } \beta' \leq 1, \\ B + X^c & \text{if } \beta' > 1. \end{cases} \end{cases}$$

LEMMA 3. *We have for all $\delta \in (0, 1]$, $p \geq 2$, $s, t \geq 0$:*

$$(59) \qquad \begin{cases} \text{(a)} \quad |B(\delta)_{t+s} - B(\delta)_t| \leq \begin{cases} Ks, & \text{if } \beta < 1, \\ Ks \log(1/\delta), & \text{if } \beta = 1, \\ Ks \delta^{1-\beta}, & \text{if } \beta > 1, \end{cases} \\ \text{(b)} \quad \mathbb{E}_t(|\widehat{X}_{t+s} - \widehat{X}_t|^p) \leq K_p s^{p/2}, \\ \text{(c)} \quad \mathbb{E}_t(|\overline{X}(\delta)^b_{t+s} - \overline{X}(\delta)^b_t|^2) \leq K s \delta^{2-\beta}, \\ \text{(d)} \quad \mathbb{E}_t(|\overline{X}(\delta)^a_{t+s} - \overline{X}(\delta)^a_t|^{\beta'}) \leq K s. \end{cases}$$

PROOF. (a) follows from (47), and (b) follows from Burkholder–Davis–Gundy inequality and Assumption 5. By (47), again, we have

$$\mathbb{E}_t(|\overline{X}(\delta)^b_{t+s} - \overline{X}(\delta)^b_t|^2) = \mathbb{E}_t\left( \int_t^{t+s} dr \int_{\{|x| \leq \delta\}} |x|^2 F'_r(dx) \right) \leq K s \delta^{2-\beta},$$

which is (c). For (d), we single out the two cases $\beta' \leq 1$ and $\beta' > 1$. In the first case, and since for any sequence $(x_m)$ we have $|\sum_m x_m|^{\beta'} \leq \sum_m |x_m|^{\beta'}$, we get

$$\mathbb{E}_t(|\overline{X}(\delta)^a_{t+s} - \overline{X}(\delta)^a_t|^{\beta'}) \leq \mathbb{E}_t\left( \left( \sum_{t < r \leq t+s} |\Delta \overline{X}(\delta)^a_r|^2 \right)^{\beta'} \right)$$

$$= \mathbb{E}_t\left( \int_t^{t+s} dr \int_{\{|x| \leq \delta\}} |x|^{\beta'} 1_\Phi(r, x) F''_r(dx) \right) \leq K s.$$

In the second case, we apply Burkholder–Davis–Gundy inequality with the exponent $1 < \beta' < 2$ to get

$$\mathbb{E}(|\overline{X}(\delta)^a_{t+s} - \overline{X}(\delta)^a_t|^{\beta'} \mid \mathcal{F}_t) \leq \mathbb{E}\left(\left(\sum_{t < r \leq t+s} |\Delta\overline{X}(\delta)^a_r|^2\right)^{\beta'/2} \Big| \mathcal{F}_t\right)$$

$$\leq \mathbb{E}\left(\sum_{t < r \leq t+s} |\Delta\overline{X}(\delta)^a_r|^{\beta'} \Big| \mathcal{F}_t\right),$$

which, exactly as before, is smaller than $Ks$, and (d) is proved. $\quad\square$

Next, we give a general result on counting processes. Let $N$ be a counting process (i.e., right continuous with $N_0 = 0$, piecewise constant, with jumps equal to 1) adapted to $(\mathcal{F}_t)$ and with predictable compensator of the form $\Lambda_t = \int_0^t \lambda \, ds$.

LEMMA 4.  *With $N$ and $\Lambda$ as above, assume further that $\lambda_t \leq u$ for some constant $u > 0$. Then, we have*

(60)   $|\mathbb{P}_t(N_{t+s} - N_t = 1) - \mathbb{E}_t(\Lambda_{t+s} - \Lambda_t)| + \mathbb{P}_t(N_{t+s} - N_t \geq 2) \leq (us)^2.$

PROOF.  Let $T_1, T_2, \ldots$ be the successive jump times of $N$ after time $(i-1)\Delta_n$. We have

$$\mathbb{P}_t(N_{t+s} - N_t \geq 1) \leq \mathbb{E}_t(N_{t+s} - N_t) = \mathbb{E}_t(\Lambda_{t+s} - \Lambda_t) \leq us,$$

$$\mathbb{P}_t(N_{t+s} - N_t = 1) = \mathbb{E}_t(N_{(t+s)\wedge T_1} - N_t)$$

$$= \mathbb{E}_t(\Lambda_{(t+s)\wedge T_1} - \Lambda_t),$$

$$\mathbb{E}_t(\Lambda_{t+s} - \Lambda_t) - \mathbb{P}_t(N_{t+s} - N_t = 1) = \mathbb{E}_t(\Lambda_{t+s} - \Lambda_{(t+s)\wedge T_1})$$

$$\leq us\mathbb{P}_t(N_{t+s} - N_t \geq 1) \leq (us)^2.$$

This gives us the first estimate. Next,

$$\mathbb{P}_t(N_{t+s} - N_t \geq 2)$$

$$= \mathbb{P}_t(T_2 \leq t + s) = \mathbb{E}_t(1_{\{T_1 < t+s\}}\mathbb{P}(T_2 \leq t + s \mid \mathcal{F}_{T_1}))$$

$$= \mathbb{E}_t\left(1_{\{T_1 < t+s\}}\mathbb{E}\left(\int_{T_1}^{(t+s)\wedge T_2} \lambda_r \, dr \Big| \mathcal{F}_{T_1}\right)\right)$$

$$\leq us\,\mathbb{P}_t(N_{t+s} - N_t \geq 1) \leq (us)^2,$$

hence the second estimate. $\quad\square$

LEMMA 5.   *With the notation* $N(\delta)_t = \sum_{s \leq t} 1_{\{|\Delta X_s| > \delta\}}$, *for all* $\delta \in (0, 1]$, $\zeta \in (0, \frac{1}{2})$ *and* $p \geq 2$ *we have*

$$(61) \qquad \begin{aligned} & \mathbb{P}_t(N(\delta)_{t+s} - N(\delta)_t \geq 1, |X(\delta)'_{t+s} - X(\delta)'_t| > \delta\zeta) \\ & \leq K_p \frac{s^{p/2}}{\zeta^p \, \delta^p} + K \frac{s^2}{\zeta^2 \, \delta^{2\beta}}. \end{aligned}$$

PROOF.   (59) and Bienaymé–Tchebycheff inequality yield

$$\mathbb{P}_t \left( |\widehat{X}_{t+s} - \widehat{X}_t| > \frac{\delta\zeta}{3} \right) \leq K_p \frac{s^{p/2}}{\zeta^p \delta^p}.$$

Next, (59) yields $|B(\delta)_{t+s} - B(\delta)_t| \leq K_0 s/\delta$ for some constant $K_0$, so we have either $|B(\delta)_{t+s} - B(\delta)_t| \leq \delta\zeta/3$ or $Ks^{p/2}/\zeta^p \delta^p \geq 1$ for some constant $K$ independent of $\zeta$ and $\delta$. Then, it remains to prove that, with the notation $M(\delta) = (x1_{\{|x| \leq \delta\}}) \star (\mu - \nu)$,

$$(62) \quad \mathbb{P}_t \left( N(\delta)_{t+s} - N(\delta)_t \geq 1, |M(\delta)_{t+s} - M(\delta)_t| > \frac{\delta\zeta}{3} \right) \leq K \frac{s^2}{\zeta^2 \, \delta^{2\beta}}.$$

For simplicity, write $N_s = N(\delta)_{t+s} - N(\delta)_t$ and $M_s = M(\delta)_{t+s} - M(\delta)_t$. By Bienaymé–Tchebycheff inequality, the left-side of (62) is not bigger than $9\mathbb{E}_t(N_s M_s^2)/\delta^2 \zeta^2$. Now, $N$ is a counting process and $M$ is a purely discontinuous square-integrable martingale, and they have no common jumps, so Itô's formula yields

$$N_s M_s^2 = 2 \int_0^s N_{r-} M_{r-} \, dM_r + \int_0^s M_{r-}^2 \, dN_r + \sum_{r \leq s} N_{r-} (\Delta M_r)^2.$$

Moreover, the compensator $N$ is as in the previous lemma, with $\lambda_s \leq K\delta^{-\beta}$, and the predictable quadratic variation of $M$ is $\Lambda'_s = \int_0^s \lambda'_r \, dr$ with $\lambda'_r \leq K\delta^{2-\beta}$ by (47). Then, taking expectations in the above display, and since the first term of the right side above is a martingale, we get

$$\begin{aligned} \mathbb{E}_t(N_s \, M_s^2) &= \mathbb{E}_t \left( \int_0^s M_r^2 \, d\Lambda_r + \int_0^s N_r \, d\Lambda'_r \right) \\ &\leq K\delta^{-\beta} \int_0^s \mathbb{E}_t(M_r^2 + \delta^2 N_r) \, dr \\ &= K\delta^{-\beta} \int_0^s \mathbb{E}_t(\Lambda'_r + \delta^2 \Lambda_r) \, dr \leq K\delta^{2(1-\beta)} s^2. \end{aligned}$$

(62) then follows.   □

LEMMA 6.   *Let* $\alpha > 0$, $\varpi \in (0, \frac{1}{2})$ *and* $\eta \in (0, \frac{1}{2} - \varpi)$, *and set*

$$(63) \qquad \rho = \eta \wedge (\varpi(\beta - \beta') - \beta'\eta) \wedge (\varpi\gamma) \wedge (1 - \varpi\beta - 2\eta).$$

*There is a constant $K$ depending on $(\alpha, \varpi, \eta)$, and also on the characteristics of $X$, such that for all $s \in (0, \Delta_n]$ and $t \geq 0$ we have*

$$(64) \qquad \left| \mathbb{P}_t(|X_{t+s} - X_t| > \alpha \Delta_n^\varpi) - \mathbb{E}_t\left( \int_t^{t+s} \overline{F}_r(\alpha \Delta_n^\varpi) \, dr \right) \right| \leq K \Delta_n^{1 - \varpi\beta + \rho},$$

$$(65) \qquad \mathbb{P}_t(\alpha \Delta_n^\varpi < |X_{t+s} - X_t| \leq \alpha \Delta_n^\varpi (1 + \Delta_n^\eta)) \leq K \Delta_n^{1 - \varpi\beta + \rho},$$

$$(66) \qquad \mathbb{P}_t(|X_{t+s} - X_t| > \alpha \Delta_n^\varpi) \leq K \Delta_n^{1 - \varpi\beta}.$$

PROOF. It is clearly enough to prove the estimates for all $\Delta_n$ small enough. We write $\delta_n = \alpha \Delta_n^\varpi$.

(1) Apply (59) and Bienaymé–Tchebycheff inequality to obtain

$$\mathbb{P}_t(|\widehat{X}_{t+s} - \widehat{X}_t| > \alpha \Delta_n^{\varpi + \eta}/3) \leq K_p \Delta_n^{p(1 - 2\varpi - 2\eta)/2},$$

$$\mathbb{P}_t(|\overline{X}(\delta_n)_{t+s}^a - \overline{X}(\delta_n)_t^a| > \alpha \Delta_n^{\varpi + \eta}/3) \leq K \Delta_n^{1 - \varpi\beta' - \eta\beta'}.$$

Moreover, for all values of $\beta$, we have $|B(\delta_n)_{t+s} - B(\delta_n)_t| \leq K \Delta_n^{-\varpi(\beta-1)^+} \times s \log(1/s)$, which is smaller than $\alpha \Delta_n^{\varpi + \eta}/3$ as soon as $\Delta_n$ is small enough because $\eta < 1 - \varpi\beta$. Since $1 - 2\varpi - 2\eta > 0$, by choosing $p$ large enough, we deduce from the previous estimates that, for all $\Delta_n$ small enough,

$$(67) \qquad \begin{aligned} &\mathbb{P}_t(|X(\delta_n)'_{t+s} + \overline{X}(\delta_n)_{t+s}^b - (X(\delta_n)'_t + \overline{X}(\delta_n)_t^b)| > \alpha \Delta_n^{\varpi + \eta}) \\ &\leq K \Delta_n^{1 - \varpi\beta + \rho}. \end{aligned}$$

(2) By Assumption 6, we have $F'_r(dx) \leq (L'/|x|^{1+\beta}) \, dx$ for some constant $L'$. We fix $n$. For each $\omega \in \Omega$, we endow the canonical (Skorokhod) space $(\Omega', \mathcal{F}', (\mathcal{F}'_t))$ of all càdlàg functions on $\mathbb{R}_+$ starting from 0 with the (unique) probability measure $Q_\omega^n$, under which the canonical process $X'$ is a semimartingale with characteristics $(0, 0, \nu'_{\delta_n, \omega})$, where

$$(68) \qquad \nu'_{\delta_n, \omega}(dr, dx) = dr \, 1_{\{|x| \leq \delta_n\}} \left( \frac{L'}{|x|^{1+\beta}} \, dx - F'_r(\omega, dx) \right).$$

This measure does not depend on $\omega'$. Hence, under $Q_\omega^n$ the process $X'$ has independent increments. $\nu'_{\delta_n, \omega}(dr, dx)$ depends measurably on $\omega$. Hence, $Q_\omega^n(d\omega')$ is a transition probability from $(\Omega, \mathcal{F})$ into $(\Omega', \mathcal{F}')$. Then, we extend $X$, $X'$ and other quantities defined on $\Omega$ or $\Omega'$ in the usual way (without changing the symbols) to the product $\widetilde{\Omega} = \Omega \times \Omega'$ endowed with the product $\sigma$-field $\widetilde{\mathcal{F}}$, the product filtration $(\widetilde{\mathcal{F}}_t)$ and the probability measure $\widetilde{\mathbb{P}}_n(d\omega, d\omega') = \mathbb{P}(d\omega) Q_\omega^n(d\omega')$.

Because of (68) and (47), and as in Lemma 3, $\mathbb{E}_{Q_\omega^n}(|X'_{t+s} - X'_t|^2 \mid \mathcal{F}'_t) \leq K s \delta_n^{2-\beta}$. Then, we see that, for some constant $C$ depending on $\alpha$ and $\beta$ but not on $n$, $i$ and $\omega$, we have

$$(69) \qquad Q_\omega^n(|X'_{t+s} - X'_t| > \alpha \Delta_n^{\varpi + \eta} \mid \mathcal{F}'_t) \leq C \Delta_n^{1 - \varpi\beta - 2\eta}.$$

(3) By well known results on extensions of spaces [see, e.g., Jacod and Shiryaev (2003), Section II.7; note that the present extension of the original space is a "very good extension"], $X'$ is a semimartingale on the extension with characteristics $(0, 0, \nu'_{\delta_n})$, where $\nu'_{\delta_n}((\omega, \omega'), dr, dx) = \nu'_{\delta_n, \omega}(dr, dx)$ and any semimartingale on the original space is a semimartingale on the extension with the same characteristics. Moreover, $X$ and $X'$ have almost surely no common jump, so the sum $Y(\delta_n)' = \overline{X}(\delta_n)^b + X'$ is a semimartingale with characteristics $(0, 0, \nu''_{\delta_n})$, where

$$\nu''_{\delta_n}(dr, dx) = dr\, 1_{\{|x| \le \delta_n\}} F'_r(dx) + \nu'_{\delta_n}(dr, dx) = 1_{\{|x| \le \delta_n\}} \frac{L'}{|x|^{1+\beta}}\, dr\, dx,$$

where the last equality comes from (68). It follows that $Y(\delta_n)'$ is a Lévy process with Lévy measure given above, or in other words it is a version of the process $Y(\delta_n)'$ of (48) with $A = L'$. Hence, we deduce from (49) and from the Lévy property of $Y(\delta_n)'$ that, for any $A \in \mathcal{F}_t$,

$$(70) \quad \widetilde{\mathbb{P}}_n(A \cap \{|Y(\delta_n)'_{t+s} - Y(\delta_n)'_t| > \alpha \Delta_n^{\varpi}(1 - 2\Delta_n^{\eta})\}) \le K\Delta_n^{4/3 - 4\varpi\beta/3} \mathbb{P}(A).$$

Next, for all $\Delta_n$ small enough, so that $C\Delta_n^{1-\varpi\beta-2\eta} \le 1/2$, we can write

$$\widetilde{\mathbb{P}}(A \cap \{|Y(\delta_n)'_{t+s} - Y(\delta_n)'_t| > \alpha \Delta_n^{\varpi}(1 - 2\Delta_n^{\eta})\})$$
$$\ge \widetilde{\mathbb{P}}(A \cap \{|\overline{X}(\delta_n)^b_{t+s} - \overline{X}(\delta_n)^b_t| > \alpha \Delta_n^{\varpi}(1 - \Delta_n^{\eta})\}$$
$$\cap \{|X'_{t+s} - X'_t| \le \alpha \Delta_n^{\varpi+\eta}\})$$
$$= \widetilde{\mathbb{E}}(1_{A \cap \{|\overline{X}(\delta_n)^b_{t+s} - \overline{X}(\delta_n)^b_t| > \alpha \Delta_n^{\varpi}(1 - \Delta_n^{\eta})\}} Q^n_{\cdot}(|X'_{t+s} - X'_t| \le \alpha \Delta_n^{\varpi+\eta}))$$
$$\ge \tfrac{1}{2} \mathbb{P}(A \cap \{|\overline{X}(\delta_n)^b_{t+s} - \overline{X}(\delta_n)^b_t| > \alpha \Delta_n^{\varpi}(1 - \Delta_n^{\eta})\}),$$

where the last inequality comes from (69). Then, by (70) and by the fact that $A$ is arbitrary in $\mathcal{F}_{(i-1)\Delta_n}$, we deduce (since necessarily $\rho \le \frac{1-\varpi\beta}{3}$) that

$$\mathbb{P}_t(|\overline{X}(\delta_n)^b_{t+s} - \overline{X}(\delta_n)^b_t| > \alpha \Delta_n^{\varpi}(1 - \Delta_n^{\eta})) \le K\Delta_n^{4/3 - 4\varpi\beta/3} \le K\Delta_n^{1-\varpi\beta+\rho}.$$

In turn, combining this with (67), we readily obtain that, for all $\Delta_n$ small enough,

$$(71) \quad \begin{aligned} &\mathbb{P}_t(|X(\delta_n)'_{t+s} - X(\delta_n)'_t| > \alpha \Delta_n^{\varpi}) \\ &\quad \le K\Delta_n^{1-\varpi\beta}(\Delta_n^{1/3 - \varpi\beta/3} + \Delta_n^{\varpi(\beta-\beta')-2\eta}) \\ &\quad \le K\Delta_n^{1-\varpi\beta+\rho}. \end{aligned}$$

(4) Now, we write $\delta'_n = \alpha \Delta_n^{\varpi}(1 + \Delta_n^{\eta})$ and also

$$\theta_{t,s} = \mathbb{E}_t\left(\int_t^{t+s} \overline{F}_r(\delta_n)\, dr\right), \qquad \theta'_{t,s} = \mathbb{E}_t\left(\int_t^{t+s} \overline{F}_r(\delta'_n)\, dr\right),$$

and the following two counting process

$$N_t^n = \sum_{s \le t} 1_{\{|\Delta X_s| > \delta_n\}}, \qquad N_t'^n = \sum_{s \le t} 1_{\{|\Delta X_s| > \delta_n'\}}.$$

Their predictable compensators are $\int_0^t \overline{F}_r(\delta_n) \, dr$ and $\int_0^t \overline{F}_r(\delta_n') \, dr$, whereas both $\overline{F}_r(\delta_n)$ and $\overline{F}_r(\delta_n')$ are smaller than $K/\Delta_n^{\varpi\beta}$. Hence, (60) gives

$$(72) \qquad \begin{cases} |\mathbb{P}_t(N_{t+s}^n - N_t^n = 1) - \theta_{t,s}| + \mathbb{P}_t(N_{t+s}^n - N_t^n \ge 2) \le K\Delta_n^{2(1-\varpi\beta)}, \\ |\mathbb{P}_t(N_{t+s}'^n - N_t'^n = 1) - \theta_{t,s}'| \le K\Delta_n^{2(1-\varpi\beta)}. \end{cases}$$

Since $N^n - N'^n$ is nondecreasing, we have

$$\mathbb{P}_t(N_{t+s}^n - N_t^n = 1, N_{t+s}'^n - N_t'^n = 0)$$
$$= \mathbb{P}_t(N_{t+s}^n - N^t = 1) - \mathbb{P}_t(N_{t+s}'^n - N_t'^n = 1)$$
$$+ \mathbb{P}_t(N_{t+s}^n - N_t^n \ge 2, N_{t+s}'^n - N_t'^n = 1).$$

Then, (72) yields

$$(73) \qquad |\mathbb{P}_t(N_{t+s}^n - N_t^n = 1, N_{t+s}'^n - N_t'^n = 0) - (\theta_{t,s} - \theta_{t,s}')| \le K\Delta_n^{2(1-\varpi\beta)}.$$

Moreover, (47) clearly implies $\theta_{t,s} - \theta_{t,s}' \le K\Delta_n^{1-\varpi\beta}(\Delta_n^\eta + \Delta_n^{\varpi(\gamma \wedge (\beta - \beta'))}) \le K\Delta_n^{1-\varpi\beta+\rho}$. We then deduce from (73) that

$$(74) \qquad \mathbb{P}_t(N_{t+s}^n - N_t^n = 1, N_{t+s}'^n - N_t'^n = 0) \le K\Delta_n^{1-\varpi\beta+\rho}.$$

(5) If $N_{t+s}^n - N_t^n = N_{t+s}'^m - N_t'^m = 1$ and $|X_{t+s} - X_t| \le \delta_n$, then, necessarily, $|X(\delta_n)_{t+s}' - X(\delta_n)_t'| > \alpha\Delta_n^{\varpi+\eta}$. Hence,

$$\mathbb{P}_t(N_{t+s}^n - N_t^n = 1, |X_{t+s} - X_t| \le \alpha\Delta_n^\varpi)$$
$$\le \mathbb{P}_t(N_{t+s}^n - N_t^n = 1, N_{t+s}'^m - N_t'^m = 0)$$
$$+ \mathbb{P}_t(N_{t+s}^n - N_t^n = 1, |X(\delta_n)_{t+s}' - X(\delta_n)_t'| > \alpha\Delta^{\varpi+\eta}).$$

Then, if we apply (61) with $p$ large enough and $\delta = \delta_n$ and $\zeta = \Delta_n^\varpi$, together with (74), we deduce that, as soon as $\Delta_n$ is small enough,

$$(75) \qquad \mathbb{P}_t(N_{t+s}^n - N_t^n = 1, |X_{t+s} - X_t| \le \alpha\Delta_n^\varpi) \le K\Delta_n^{1-\varpi\beta+\rho}.$$

Finally, $X_{t+s} - X_t = X(\delta_n)_{t+s}' - X(\delta_n)_t'$ on the set $\{N_{t+s}^n - N_t^n = 0\}$, so
$$\mathbb{P}_t(|X_{t+s} - X_t| > \alpha\Delta_n^\varpi) = \mathbb{P}_t(N_{t+s}^n - N_t^n = 1)$$
$$- \mathbb{P}_t(N_{t+s}^n - N_t^n = 1, |X_{t+s} - X_t| \le \alpha\Delta_n^\varpi)$$
$$+ \mathbb{P}_t(N_{t+s}^n - N_t^n = 0, |X(\delta_n)_{t+s}' - X(\delta_n)_t'| > \alpha\Delta_n^\varpi)$$
$$+ \mathbb{P}_{i-1}^n(N_{t+s}^n - N_t^n \ge 2, |X_{t+s} - X_t| > \alpha\Delta_n^\varpi).$$

Then, if we combine (71), (72) and (75), we get that for all $\Delta_n$ small enough, we readily obtain (64). We also trivially deduce (66) from (47) and (64).

(6) Finally, a close look at the previous argument shows that (64) also holds with $\alpha\Delta_n^\varpi(1+\Delta_n^\eta)$ and $\theta_{t,s}i'$ in place of $\alpha\Delta_n^\varpi$ and $\theta_{t,s}$. Therefore, (65) follows upon using the property $\theta_{t,s} - \theta'_{t,s} \leq K\Delta_n^{1-\varpi\beta+\rho}$ proved above.  □

LEMMA 7.  *Under the assumption and with the notation of Lemma 6, and if $M$ is a bounded martingale, we have (with $K$ depending also on $M$, recall $s \leq \Delta_n$)*

$$
\begin{aligned}
&|\mathbb{E}_t((M_{t+s} - M_t)1_{\{|X_{t+s}-X_t|>\alpha\Delta_n^\varpi\}})| \\
(76) \qquad &\leq K\Delta_n^{1-\varpi\beta+\rho} + K\Delta_n^{1-(\varpi+\eta)\beta}\mathbb{E}_t(|M_{t+s} - M_t|) \\
&\quad + K\Delta_n^{(1-(\varpi+\eta)\beta)/2}\sqrt{\mathbb{E}_t(|M_{t+s} - M_t|^2)}.
\end{aligned}
$$

PROOF.  (1) There exist $C^2$ functions $f_n$ such that

$$
(77) \qquad \begin{cases} 1_{\{|x|>\alpha\Delta_n^\varpi(1+2\Delta_n^\eta/3)\}} \leq f_n(x) \leq 1_{\{|x|>\alpha\Delta_n^\varpi(1+\Delta_n^\eta/3)\}} \\ |f'_n(x)| \leq \dfrac{K}{\Delta_n^{\varpi+\eta}}, \qquad |f''_n(x)| \leq \dfrac{K}{\Delta_n^{2(\varpi+\eta)}}. \end{cases}
$$

With $\widehat{X}' = X - B - X^c$, and since $M$ is bounded, we have

$$
\begin{aligned}
&|\mathbb{P}_t((M_{t+s} - M_t)1_{\{|X_{t+s}-X_t|>\alpha\Delta_n^\varpi\}}) - \mathbb{E}_t((M_{t+s} - M_t)\,f_n(\widehat{X}'_{t+s} - \widehat{X}'_t))| \\
(78) \qquad &\leq K\mathbb{P}_t(\alpha\Delta_n^\varpi < |X_{t+s} - X_t| \leq \alpha\Delta_n^\varpi(1+\Delta_n^\eta)) \\
&\quad + K\mathbb{E}_t(|f_n(X_{t+s} - X_t) - f_n(\widehat{X}'_{t+s} - \widehat{X}'_t)|).
\end{aligned}
$$

Now, we have

$$
|f_n(x+y) - f_n(x)| \leq 1_{\{|y|>\alpha\Delta_n^{\varpi+\eta}/3\}} + \frac{K}{\Delta_n^{\varpi+\eta}}|y|1_{\{\alpha\Delta_n^\varpi < |x+y| \leq \alpha\Delta_n^\varpi(1+\Delta_n^\eta)\}}.
$$

If we apply this with $x = \widehat{X}'_{t+s} - \widehat{X}'_t$ and $y = (B+X^c)_{t+s} - (B+X^c)_t$, plus (59)(b) for $p$ large enough and Bienaymé–Tchebycheff inequality and $1 - 2\varpi - 2\eta > 0$, plus (65) and (59)(b) again and Hölder's inequality, we obtain that the right side of (78) is smaller than $K\Delta_n^{1-\varpi\beta+\rho}$. Therefore, it remains to prove that

$$
\begin{aligned}
&|\mathbb{E}_t((M_{t+s} - M_t)\,f_n(\widehat{X}'_{t+s} - \widehat{X}'_t))| \\
(79) \qquad &\leq K\Delta_n^{1-\varpi\beta+\rho} + K\Delta_n^{1-(\varpi+\eta)\beta}\mathbb{E}_t(|M_{t+s} - M_t|) \\
&\quad + K\Delta_n^{(1-(\varpi+\eta)\beta)/2}\sqrt{\mathbb{E}_t(|M_{t+s} - M_t|^2)}.
\end{aligned}
$$

(2) According to Theorem III.4.20 of Jacod and Shiryaev (2003), we can "project" the martingale $M$ onto the random measure $\mu$, which amounts to decomposing it as the sum of two martingales $M = M' + M''$ with the following properties:

$$(80) \quad \begin{cases} \bullet \ M' = \delta \star (\mu - \nu) \text{ for some predictable function } \delta, \\ \bullet \ \sum_{s \le t} \phi(s, \Delta X_s) \Delta M''_s \text{ is a martingale as soon as} \\ \quad \phi \text{ is a predictable function satisfying } |\phi(\omega, s, x)| \le K(1 \wedge |x|), \\ \bullet \ \text{the difference of the predictable quadratic variations} \\ \quad \langle M, M \rangle - \langle M', M' \rangle \text{ is nondecreasing.} \end{cases}$$

Note that we may even choose $\delta$ bounded, because $M$ is bounded [the fact that all the above processes are martingales and not only local martingales comes from the boundedness of $M$ and of $\int (x^2 \wedge 1) F_t(dx)$]. We will also use the following consequence of the third property above:

$$(81) \quad \begin{aligned} \mathbb{E}_t &\left( \int_t^{t+s} du \int F_{t+u}(x) \delta(t+u, x)^2 \right) \\ &= \mathbb{E}_t(\langle M', M' \rangle_{t+s} - \langle M', M' \rangle_t) \\ &\le \mathbb{E}_t((M_{t+s} - M_t)^2). \end{aligned}$$

With $t$ being fixed below, for simplicity, we write $Y_r = \widehat{X}'_{t+r} - \widehat{X}'_t$. Since $M$ is a bounded martingale and $Y$ a semimartingale with vanishing continuous martingale part, and $f_n(Y)$ is bounded, we deduce, from Itô's formula and the properties (80), that the product $(M_{t+r} - M_t) f_n(Y_r)$ is the sum of a martingale plus the process $\int_0^r \gamma_u^n \, du$, where

$$\gamma_u^n = \int F_{t+u}(dx)((M_{t+u} - M_t) g_n(Y_u, x) + \delta(t+u, x) h_n(Y_u, x)),$$

where

$$h_n(y, x) = f_n(y + x) - f_n(y), \qquad g_n(y, x) = h_n(x, y) - f'_n(y) x \mathbf{1}_{\{|x| \le 1\}}.$$

An easy computation allows us to deduce of (77) that

$$|h_n(y, x)| \le K \frac{|x| \wedge 1}{\Delta_n^{\varpi + \eta}},$$

$$\begin{aligned} |g_n(y, x)| \le \ &\mathbf{1}_{\{|x| > \alpha \Delta_n^{\varpi + \eta}\}} + K \mathbf{1}_{\{\alpha \Delta_n^{\varpi} < |y| \le \alpha \Delta_n^{\varpi}(1 + \Delta_n^n)\}} \\ &\times \left( \frac{x^2}{\Delta_n^{2\varpi + 2\eta}} \mathbf{1}_{\{|x| \le \alpha \Delta_n^{\varpi + \eta}\}} + \frac{|x| \wedge 1}{\Delta_n^{\varpi + \eta}} \mathbf{1}_{\{|x| > \alpha \Delta_n^{\varpi + \eta}\}} \right). \end{aligned}$$

Now, we apply the first estimate of (47) with $x = \alpha \Delta^{\varpi + \eta}$ and the second and third ones with $u = \alpha \Delta^{\varpi + \eta}$ plus Cauchy–Schwarz inequality, to get, for any $\varepsilon > 0$ (recall that $\delta$ is bounded),

$$\begin{aligned} |\gamma_u^n| \le \ &K \Delta_n^{-(\varpi + \eta)\beta} |M_{t+u} - M_t| + K \Delta_n^{-(\varpi + \eta)\beta/2} \left( \int F_{t+u}(x) \delta(t+u, x)^2 \right)^{1/2} \\ &+ K_\varepsilon |M_{t+u} - M_t| \, \Delta_n^{-(\beta + \varepsilon)(\varpi + \eta)} \mathbf{1}_{\{\alpha \Delta_n^{\varpi} < |Y_u| \le \alpha \Delta_n^{\varpi}(1 + \Delta_n \eta)\}}. \end{aligned}$$

Since $\eta < 1/2 - \varpi$, we have $\beta(\varpi + \eta) < 1$ and, thus, $(\beta + \varepsilon)(\varpi + \eta) = 1$ for a suitable $\varepsilon > 0$. Moreover, $\mathbb{E}_t(|Z_u|) \le \mathbb{E}_t(|Z_s|)$ if $u \le s$, because $Z$ is a martingale. Therefore, since $M$ is bounded and $s \le \Delta_n$, we get

$$
\begin{aligned}
&|\mathbb{E}_t((M_{t+s} - M_t) f_n(\widehat{X}'_{t+s} - \widehat{X}'_t))| \\
&\quad = \left| \mathbb{E}_t\left( \int_0^s \gamma_r^n \, dr \right) \right| \le \int_0^s \mathbb{E}_t(|\gamma_r^n|) \, dr \\
&\quad \le K\Delta_n^{1-(\varpi+\eta)\beta} \mathbb{E}_t(|M_{t+s} - M_t|) \\
&\qquad + K\Delta_n^{(1-(\varpi+\eta)\beta)/2} \left( \mathbb{E}_t\left( \int_t^{t+s} du \int F_{t+u}(x)\delta(t+u,x)^2 \right) \right)^{1/2} \\
&\qquad + K\Delta_n^{-1} \int_0^{\Delta_n} \mathbb{P}_t(\alpha\Delta_n^\varpi < |Y_r| \le \alpha\Delta_n^\varpi(1 + \Delta_n\eta)) \, dr,
\end{aligned}
$$

where we have used Cauchy–Schwarz inequality. By (65), for the process $\widehat{X}'$ instead of $X$ and by (81), we readily deduce (79). $\quad\square$

8.3. *Some auxiliary limit theorems.* Below, recall the process $\overline{A}$ of (11). We still assume Assumptions 5 and 6 and also $\varpi \in (0, \frac{1}{2})$ and $\alpha > 0$.

LEMMA 8. *Let* $\rho' = \frac{1}{2} \wedge (\varpi(\beta - (\beta - \gamma) \vee \beta'))$. *Then, for all* $t > 0$, *the sequence*

$$
(82) \qquad \left( \Delta_n^{-\rho'} \left| \sum_{i=1}^{[t/\Delta_n]} \Delta_n^{\varpi\beta} \mathbb{E}_{i-1}^n \left( \int_{(i-1)\Delta_n}^{i\Delta_n} \overline{F}_t(\alpha\Delta_n^\varpi) \, dt \right) - \frac{\overline{A}_t}{\alpha^\beta} \right| \right)_{n \ge 1}
$$

*is tight.*

PROOF. Let $\theta_i^n = \int_{(i-1)\Delta_n}^{i\Delta_n} \overline{F}_t(\alpha\Delta_n^\varpi) \, dt$ and $\eta_i^n = \int_{(i-1)\Delta_n}^{i\Delta_n} A_t \, dt$. We deduce from (47) that

$$
\left| \Delta_n^{\varpi\beta} \theta_i^n - \frac{1}{\alpha^\beta} \eta_i^n \right| \le K\Delta_n^{1+\varpi(\beta-(\beta-\gamma)\vee\beta')}.
$$

Then, obviously,

$$
\mathbb{E}\left( \Delta_n^{-\rho'} \sum_{i=1}^{[t/\Delta_n]} \mathbb{E}_{i-1}^n \left( \left| \Delta_n^{\varpi\beta} \theta_i^n - \frac{1}{\alpha^\beta} \eta_i^n \right| \right) \right) \le Kt,
$$

and, since $A_t$ is bounded, we have $|\bar{A}_t - \sum_{i=1}^{[t/\Delta_n]} \eta_i^n| \le Kt\Delta_n$, whereas $\rho' < 1$. It thus remains to prove that

$$
(83) \qquad \text{the sequence } \left( \Delta_n^{-\rho'} \left| \sum_{i=1}^{[t/\Delta_n]} (\eta_i^n - \mathbb{E}_{i-1}^n(\eta_i^n)) \right| \right)_{n \ge 1} \text{ is tight.}
$$

Since $\zeta_i^n = \Delta_n^{-\rho'}(\eta_i^n - \mathbb{E}_{i-1}^n(\eta_i^n))$ is a martingale increment, for (83), it is enough to check that $a_n(t) = \mathbb{E}(\sum_{i=1}^{[t/\Delta_n]}(\zeta_i^n)^2)$ is bounded. However, since $A_t$ is bounded, we have $|\zeta_i^n|^2 \le K\Delta_n^{2-2\rho'}$, so $a_n(t) \le Kt\Delta_n^{1-2\rho'} \le K$ because $\rho' \le 1/2$. $\square$

LEMMA 9. (a) *Let*

$$\chi' = (\varpi\gamma) \wedge \frac{1-\varpi\beta}{3} \wedge \frac{\varpi(\beta-\beta')}{1+\beta'} \wedge \frac{1-2\varpi}{2}. \tag{84}$$

*Then, for all $\varepsilon > 0$ and all $t > 0$ the sequence of variables*

$$\left(\Delta_n^{\varepsilon-\chi'}\left|\Delta_n^{\varpi\beta}\sum_{i=1}^{[t/\Delta_n]}\mathbb{P}_{i-1}^n(|\Delta_i^n X| > \alpha\Delta_n^{\varpi}) - \frac{\overline{A_t}}{\alpha^\beta}\right|\right)_{n\ge1} \tag{85}$$

*is tight, and, in particular, we have*

$$\Delta_n^{\varpi\beta}\sum_{i=1}^{[t/\Delta_n]}\mathbb{P}_{i-1}^n(|\Delta_i^n X| > \alpha\Delta_n^{\varpi}) \xrightarrow{\mathbb{P}} \frac{\overline{A_t}}{\alpha^\beta}. \tag{86}$$

(b) *If further $\beta' < \frac{\beta}{2+\beta}$ and $\gamma > \frac{\beta}{2}$ and $\varpi < \frac{1}{2+\beta} \wedge \frac{2}{5\beta}$, and if $M$ is a bounded continuous martingale, we also have*

$$\Delta_n^{-\varpi\beta/2}\left|\Delta_n^{\varpi\beta}\sum_{i=1}^{[s/\Delta_n]}\mathbb{P}_{i-1}^n(|\Delta_i^n X| > \alpha\Delta_n^{\varpi}) - \frac{\overline{A_s}}{\alpha^\beta}\right| \xrightarrow{\mathbb{P}} 0, \tag{87}$$

$$\Delta_n^{\varpi\beta/2}\sum_{i=1}^{[t/\Delta_n]}|\mathbb{E}_{i-1}^n(\Delta_i^n M 1_{\{|\Delta_i^n X| > \alpha\Delta_n^{\varpi}\}})| \xrightarrow{\mathbb{P}} 0. \tag{88}$$

PROOF. (a) Let $\eta \in (0, \frac{1}{2} - \varpi)$ and $\rho$ given by (63) and $\rho'$ like in Lemma 8. From (64) and (82) we deduce the tightness of the sequence (85), provided we substitute $\chi' - \varepsilon$ with $\rho' \wedge \rho$. It is thus enough to show that we can choose $\eta \le \frac{1-2\varpi}{2}$ such that $\rho' \wedge \rho = \chi' - \varepsilon$, and this is achieved by taking $\eta = \frac{1-\varpi\beta}{3} \wedge \frac{\varpi(\beta-\beta')}{1+\beta'} \wedge \frac{1-2\varpi}{2} - \varepsilon$, as a simple computation shows.

(b) In view of (a), (87) follows from the property $\chi' > \varpi\beta/2$, which is an easy consequence of the assumptions.

It remains to prove (88). By (76), the left-hand side of this expression is smaller than

$$Kt\Delta_n^{\rho-\varpi\beta/2} + K\Delta_n^{1-\eta\beta-\varpi\beta/2}\sum_{i=1}^{[t/\Delta_n]+1}\mathbb{E}_{i-1}^n(|\Delta_i^n M|)$$

$$+ K\Delta_n^{(1-\eta\beta)/2}\sum_{i=1}^{[t/\Delta_n]+1}\sqrt{\mathbb{E}_{i-1}^n(|\Delta_i^n M|)}.$$

By the Cauchy–Schwarz inequality, this is smaller than

$$K(t + \sqrt{t})\left(\Delta_n^{\rho - \varpi\beta/2} + \Delta_n^{1/2 - \eta\beta - \varpi\beta/2}\left(\sum_{i=1}^{[t/\Delta_n]+1} \mathbb{E}_{i-1}^n(|\Delta_i^n M|^2)\right)^{1/2}\right).$$

A well known property of martingales yields

$$\mathbb{E}\left(\sum_{i=1}^{[t/\Delta_n]+1} \mathbb{E}_{i-1}^n(|\Delta_i^n M|^2)\right) = \mathbb{E}((M_{\Delta_n([t/\Delta_n]+1)} - M_0)^2),$$

which is bounded (in $n$). Therefore, we deduce that (88) holds, provided we have $\rho > \varpi\beta/2$ and also $1 - 2\eta\beta > \varpi\beta$. The first condition has already been checked, and the second is satisfied because $\eta \leq (1 - \varpi\beta)/3$. This ends the proof. $\square$

8.4. *Law of large numbers and central limit theorems.* Here, again, we assume Assumptions 5 and 6, and we have $\alpha > 0$ and $\varpi \in (0, 1/2)$.

PROPOSITION 1. *If $\chi$ is given by (19), then for each $t > 0$ and each $\varepsilon > 0$ the sequence*

$$\left(\Delta_n^{\varepsilon - \chi}\left|\Delta_n^{\varpi\beta} U(\varpi, \alpha)_t^n - \frac{\bar{A}_t}{\alpha^\beta}\right|\right)_{n \geq 1} \tag{89}$$

*is tight, and, in particular,*

$$\Delta_n^{\varpi\beta} U(\varpi, \alpha)_t^n \xrightarrow{\mathbb{P}} \frac{\bar{A}_t}{\alpha^\beta}. \tag{90}$$

PROOF. Set

$$\zeta_i^n = \Delta_n^{\varpi\beta/2}(1_{\{|\Delta_i^n X| > \alpha\Delta_n^\varpi\}} - \mathbb{P}_{i-1}^n(|\Delta_i^n X| > \alpha\Delta_n^\varpi)). \tag{91}$$

By virtue of Lemma 9, and since $\chi = \chi' \wedge (\varpi\beta/2)$, it suffices to prove that the sequence $\sum_{i=1}^{[t/\Delta_n]} \zeta_i^n$ is tight. Since the $\zeta_i^n$'s are martingale increments, it is enough to show that the sequence $a_n(t) = \sum_{i=1}^{[t/\Delta_n]} \mathbb{E}((\zeta_i^n)^2)$ is bounded. But (66) yields $\mathbb{E}((\zeta_i^n)^2) \leq K\Delta_n$, which in turn yields $a_n(t) \leq Kt$. $\square$

PROPOSITION 2. *Let $\alpha' > \alpha$. If we have $\beta' < \frac{\beta}{2+\beta}$ and $\gamma > \frac{\beta}{2}$ and $\varpi < \frac{1}{2+\beta} \wedge \frac{2}{5\beta}$, the pair of processes*

$$\Delta_n^{-\varpi\beta/2}\left(\Delta_n^{\varpi\beta} U(\varpi, \alpha)_t^n - \frac{\bar{A}_t}{\alpha^\beta}, \Delta_n^{\varpi\beta} U(\varpi, \alpha')_t^n - \frac{\bar{A}_t}{\alpha'^\beta}\right) \tag{92}$$

*converges stably in law to a continuous Gaussian martingale $(\overline{W}, \overline{W}')$ independent of $\mathcal{F}$, with*

$$\mathbb{E}(\overline{W}_t^2) = \frac{\bar{A}_t}{\alpha^\beta}, \qquad \mathbb{E}(\overline{W}'^2_t) = \frac{\bar{A}_t}{\alpha'^\beta}, \qquad \mathbb{E}(\overline{W}_t\overline{W}'_t) = \frac{\bar{A}_t}{\alpha'^\beta}. \tag{93}$$

PROOF. Define $\zeta_i^n$ by (91) and associate $\zeta_i'^n$ with $\alpha'$ in the same way. The variables $\zeta_i^n$ and $\zeta_i'^n$ are martingale increments and smaller than $K\Delta_n^{\varpi\beta/2}$. So, in view of (87), it is enough, by using a criterion for stable convergence of triangular arrays found in Jacod and Shiryaev (2003) (see Theorem IX.7.28) to prove the following:

$$(94) \qquad \begin{cases} \displaystyle\sum_{i=1}^{[t/\Delta_n]} \mathbb{E}_{i-1}^n((\zeta_i^n)^2) \overset{\mathbb{P}}{\longrightarrow} \frac{\bar{A}_t}{\alpha^\beta}, \\[2.5ex] \displaystyle\sum_{i=1}^{[t/\Delta_n]} \mathbb{E}_{i-1}^n((\zeta_i'^n)^2) \overset{\mathbb{P}}{\longrightarrow} \frac{\bar{A}_t}{\alpha'^\beta}, \\[2.5ex] \displaystyle\sum_{i=1}^{[t/\Delta_n]} \mathbb{E}_{i-1}^n(\zeta_i^n\,\zeta_i'^n) \overset{\mathbb{P}}{\longrightarrow} \frac{\bar{A}_t}{\alpha'^\beta}, \end{cases}$$

$$(95) \qquad \begin{cases} \displaystyle\sum_{i=1}^{[t/\Delta_n]} \mathbb{E}_{i-1}^n(\zeta_i^n\,\Delta_i^n M) \overset{\mathbb{P}}{\longrightarrow} 0, \\[2.5ex] \displaystyle\sum_{i=1}^{[t/\Delta_n]} \mathbb{E}_{i-1}^n(\zeta_i'^n\,\Delta_i^n M) \overset{\mathbb{P}}{\longrightarrow} 0, \end{cases}$$

where $M$ is any bounded martingale.

Recalling $\alpha < \alpha'$ and also $\mathbb{P}_{i-1}^n(|\Delta_i^n X| > \alpha\Delta_n^\varpi) \leq K\Delta_n^{1-\varpi\beta}$ by (66), we deduce (94) from (86). As for (95), and since $\mathbb{E}_{i-1}^n(\Delta_i^n M) = 0$, it readily follows from (88). □

### 8.5. *Proofs of the main theorems.*

PROOFS OF THEOREMS 1 AND 2. As said before, we can assume Assumptions 5 and 6. We will write $\eta = \chi - \varepsilon$ if we want to prove the first theorem, in which case we choose $\varepsilon$ small enough to have $\eta > 0$, and $\eta = \varpi\beta/2$ if we want to prove the second one. Then, we set

$$V_n = \Delta_n^{-\eta}(\Delta_n^{\varpi\beta} U(\varpi,\alpha)_t^n - \bar{A}_t/\alpha^\beta), \qquad V_n' = \Delta_n^{-\eta}(\Delta_n^{\varpi\beta} U(\varpi,\alpha')_t^n - \bar{A}_t/\alpha'^\beta),$$

so that on the set $\{\bar{A}_t > 0\}$ we have

$$(96) \qquad \Delta_n^{-\eta}(\hat{\beta}_n(t,\varpi,\alpha,\alpha') - \beta) = \frac{\Delta_n^{-\eta}}{\log(\alpha'/\alpha)} \log \frac{1 + \alpha^\beta\Delta_n^\eta V_n/\bar{A}_t}{1 + \alpha'^\beta\Delta_n^\eta V_n'/\bar{A}_t}.$$

In all cases, the two sequences $(V_n)$ and $(V_n')$ are tight, so the right-hand side is equivalent (in probability) to the variable $(\alpha^\beta V_n - \alpha'^\beta V_n')/(\bar{A}_t \log(\alpha'/\alpha))$ on the set $\{\bar{A}_t > 0\}$. At this stage, Theorems 1 and 2 readily follow from Lemma 9(a) and Proposition 2, respectively. □

PROOF OF THEOREM 3.   This is a trivial consequence of the properties of the stable convergence in law and of Theorem 2.  □

PROOFS OF THEOREMS 4 AND 5.   In the situation of these two theorems, we have $\bar{A}_t = At > 0$, so we can apply Theorems 1 and 3 with $\{\bar{A}_t > 0\} = \Omega$. Indeed, we have Assumptions 1 and 2 with $\beta' = 0$ and $\gamma$ arbitrary large. So, the only differences are on the conditions on $\varpi$, namely $\varpi < 1/\beta$ for the consistency and $\varpi < 2/(3\beta)$ for the CLT in the first theorem, and $\varpi < 1/(2 + \beta)$ for the CLT in the second one, instead of $\varpi < 1/2$ in Theorem 1 and $\varpi < \frac{1}{2+\beta} \wedge \frac{2}{5\beta}$ in Theorem 3.

In fact, the improvement lies in Lemma 6. When $X = Y$ we have the trivial estimate, coming from (24) and from the scaling property of $Y$,

$$(97) \qquad \mathbb{P}_t(|X_{t+s} - X_t| > x) = \mathbb{P}(|X_s| > x) = \frac{As}{x^\beta} + O(s^2/x^{2\beta}),$$

regardless of the value of $x > 0$. So, we have all claims of Lemma 6 with any $\eta \in (0, 1 - \varpi\beta]$ and with $\rho = 1 - \varpi\beta$, and, thus, (a) of Lemma 9 holds for some $\chi > 0$ as soon as $\varpi < 1/\beta$, whereas (b) of that lemma holds if further $\rho > \varpi\beta/2$, that is, when $\varpi < 2/(3\beta)$.

When $X_t = \sigma W_t + Y_t$, we obtain for any $\eta \in (0, 1/2 - \varpi)$ and $s \leq \Delta_n$

$$\mathbb{P}_t(|X_{t+s} - X_t| > \alpha\Delta_n^\varpi) = \mathbb{P}(|X_s| > \alpha\Delta_n^\varpi) = \frac{A\Delta_n^{1-\varpi\beta}}{\alpha^\beta} + O(\Delta_n^{1-\varpi\beta+\eta}),$$

[by applying the estimate (97) for $Y$ and the fact that $\mathbb{E}(|X_s^c|^p) = K_p s^{p/2}$ and Bienaymé–Tchebycheff inequality]. Then, we have all claims of Lemma 6 with any $\eta \in (0, (\frac{1}{2} - \varpi) \wedge (1 - \varpi\beta))$ and with $\rho = \eta$, and, thus, (b) of Lemma 9 holds as soon as we can find $\eta$ as above, with $\rho = \eta > \frac{\varpi\beta}{2}$. This is possible, provided $\varpi < \frac{1}{2+\beta} \wedge \frac{2}{3\beta}$. The rest of the proofs of Theorems 1 and 3 goes through, and we are finished.  □

**Acknowledgements.**   We are very grateful to two referees and an Associate Editor for many helpful comments.

Department of Economics
Princeton University and NBER
Princeton, NJ 08544-1021
USA
E-mail: yacine@princeton.edu

Institut de Mathématiques de Jussieu
CNRS UMR 7586
Université P. et M. Curie (Paris-6)
75252 Paris Cédex 05
France
E-mail: jj@ccr.jussieu.fr